\documentclass[leqno,a4paper,12pt]{article}
\usepackage{amsmath,amsthm,bbm}
\usepackage[sans]{dsfont}

\newtheorem{Thm}{\indent Theorem}[section]
\newtheorem{Prop}[Thm]{\indent Proposition}

\newtheorem{Cor}[Thm]{\indent Corollary}

\theoremstyle{definition}
\newtheorem{Def}[Thm]{\indent Definition}
\newtheorem{Rem}[Thm]{\indent Remark}
\newtheorem{Ex}[Thm]{\indent Example}

\newtheorem{Cons}[Thm]{\indent Construction}
\newtheorem{Ass}[Thm]{\indent Assumption}
\newtheorem{Ques}[Thm]{\indent Question}

\def\qed{{\hskip0pt\unskip\unskip\nobreak\hfil\penalty50
          \hskip1em\hbox{}\nobreak\hfil
          {\bf q.e.d.}%
          \parfillskip=0pt\finalhyphendemerits=0
          \par}\medskip}

\newenvironment{Proof}
               {{\it Proof.}\quad}
               {\qed}

\newenvironment{Proofof}[1]
               {{\it Proof of #1.}\quad}
               {\qed}

\newcommand{\Prime}{\kern3\fontdimen1\font$'$\kern-7\fontdimen1\font}

\long\def\forget#1{}

\long\def\beginSIDEREMARK#1\endSIDEREMARK
    {{\par\bigskip\advance\leftskip by 2cm
                  \advance\rightskip by -2cm\noindent
      {\bf Our own side remark:} #1
      \par\bigskip\noindent}}

\long\def\beginFORGET#1\endFORGET{#1}
\long\def\beginFORGET#1\endFORGET{}

\def\?{\ ???\ \immediate\write16{}%
\immediate\write16{Warning: There was still a question mark . . . }%
\immediate\write16{}}

\usepackage{amsmath}

\usepackage{exscale}
\usepackage{amssymb}

\input cyracc.def
\font\tencyr=wncyr6
\def\cyr{\tencyr\cyracc}
\newcommand{\cyrb}{{\cyr B}}

\newcommand{\B}{{\rm{B}}}
\newcommand{\C}{{\rm{C}}}

\newcommand{\BA}{{\mathbb{A}}}
\newcommand{\BB}{{\mathbb{B}}}
\newcommand{\BC}{{\mathbb{C}}}
\newcommand{\BD}{{\mathbb{D}}}

\newcommand{\BG}{{\mathbb{G}}}

\newcommand{\BQ}{{\mathbb{Q}}}
\newcommand{\BR}{{\mathbb{R}}}

\newcommand{\BZ}{{\mathbb{Z}}}

\newcommand{\Fc}{{\mathfrak{c}}}

\newcommand{\Fg}{{\mathfrak{g}}}

\newcommand{\Fn}{{\mathfrak{n}}}

\newcommand{\FA}{{\mathfrak{A}}}

\newcommand{\FS}{{\mathfrak{S}}}

\newcommand{\FZ}{{\mathfrak{Z}}}

\newcommand{\CC}{{\cal C}}
\newcommand{\CD}{{\cal D}}

\newcommand{\CH}{{\cal H}}

\newcommand{\CV}{{\cal V}}

\newcommand{\ch}{\mathop{\rm CH}\nolimits}
\newcommand{\diag}{\mathop{\rm diag}\nolimits}
\newcommand{\Spec}{\mathop{{\bf Spec}}\nolimits}

\newcommand{\End}{\mathop{\rm End}\nolimits}

\newcommand{\GL}{\mathop{\rm GL}\nolimits}
\newcommand{\Gm}{\mathop{\BG_{m,\BQ}}\nolimits}
\newcommand{\GFm}{\mathop{\BG_{m,F}}\nolimits}
\newcommand{\GLm}{\mathop{\BG_{m,L}}\nolimits}
\newcommand{\Gr}{\mathop{\rm Gr}\nolimits}

\newcommand{\Hom}{\mathop{\rm Hom}\nolimits}

\newcommand{\Rep}{\mathop{\rm Rep}\nolimits}

\newcommand{\Sym}{\mathop{\rm Sym}\nolimits}

\newcommand{\loccit}{[loc.$\;$cit.]}

\def\halb{\frac{1}{2}}

\def\id{{\rm id}}

\newbox\mybox
\def\arrover#1{\mathrel{
       \setbox\mybox=\hbox spread 1.4em{\hfil$\scriptstyle#1$\hfil}
       \vbox{\offinterlineskip\copy\mybox
             \hbox to\wd\mybox{\rightarrowfill}}}}
\def\larrover#1{\mathrel{
       \setbox\mybox=\hbox spread 1.4em{\hfil$\scriptstyle#1$\hfil}
       \vbox{\offinterlineskip\copy\mybox
             \hbox to\wd\mybox{\leftarrowfill}}}}

\def\ontoover#1{\mathrel{
       \setbox\mybox=\hbox spread 1.4em{\hfil$\scriptstyle#1$\hfil}
       \vbox{\offinterlineskip\copy\mybox
             \hbox to\wd\mybox{\rightarrowfill\hskip-2.8mm
                               $\rightarrow$}}}}
\def\leftontoover#1{\mathrel{
       \setbox\mybox=\hbox spread 1.4em{\hfil$\scriptstyle#1$\hfil}
       \vbox{\offinterlineskip\copy\mybox
             \hbox to\wd\mybox{$\leftarrow$\hskip-2.8mm
                               \leftarrowfill}}}}
\def\longto{\longrightarrow}
\def\into{\hookrightarrow}
\def\onto{\ontoover{\ }}

\def\isoto{\arrover{\sim}}

\def\longinto{\lhook\joinrel\longrightarrow}

\usepackage[curve,matrix,arrow,cmtip]{xy}
\NoComputerModernTips

\def\myxymessage{\def\messagetext
   {Here an xy-pic diagram was omitted to speed up compilation . . . }
   \immediate\write16{\messagetext}
   \hbox{\bf \messagetext}}
\def\filxymatrix#1{\myxymessage}
\def\filxyarray#1{\myxymessage}

\newdir^{ (}{{}*!/-3pt/\dir^{(}}
\newdir_{ (}{{}*!/-3pt/\dir_{(}}
\newdir^{ )}{{}*!/+3pt/\dir^{)}}
\newdir_{ )}{{}*!/+3pt/\dir_{)}}

\def\rscript#1{\hbox to 0pt{$\scriptstyle#1$\hss}}

\let\oldbullet\bullet
\def\bullet{{\mathchoice{\oldbullet}%
                        {\oldbullet}%
                        {\scriptscriptstyle\oldbullet}%
                        {\oldbullet}}}

\newcommand{\argdot}{{\;\bullet\;}}

\newcommand{\ua}{\mathop{\underline{\alpha}}\nolimits}

\newcommand{\ujast}{\mathop{j_{!*}}\nolimits}

\newcommand{\CHMM}{\mathop{CHM(M^K)}\nolimits}
\newcommand{\CHsMM}{\mathop{CHM^s(M^K)}\nolimits}

\newcommand{\CHUM}{\mathop{CHM(U)}\nolimits}
\newcommand{\CHXM}{\mathop{CHM(X)}\nolimits}

\newcommand{\CHFDM}{\mathop{CHM(\BA)_F}\nolimits}
\newcommand{\CHFMM}{\mathop{CHM(M^K)_F}\nolimits}
\newcommand{\CHQMM}{\mathop{CHM(M^K)_\BQ}\nolimits}

\newcommand{\CHFUM}{\mathop{CHM(U)_F}\nolimits}
\newcommand{\CHFWM}{\mathop{CHM(W)_F}\nolimits}
\newcommand{\CHFXM}{\mathop{CHM(X)_F}\nolimits}

\newcommand{\CHFZoM}{\mathop{CHM(W')_F}\nolimits}
\newcommand{\DBcM}{\mathop{DM_{\text{\cyrb},c}}\nolimits}
\newcommand{\DBPAbcM}{\mathop{DM_{\text{\cyrb},c,\Phi}^{Ab}}\nolimits}
\newcommand{\DBcFpPAbM}{\mathop{\DBPAbcM(\partial (M^K)^*)_F}\nolimits}
\newcommand{\DBcFZPAbM}{\mathop{\DBPAbcM(Z)_F}\nolimits}
\newcommand{\DBcFbM}{\mathop{\DBcM(\bullet)_F}\nolimits}
\newcommand{\DBcDM}{\mathop{\DBcM(\BA)}\nolimits}
\newcommand{\DBcFDM}{\mathop{\DBcM(\BA)_F}\nolimits}
\newcommand{\DBcFBsM}{\mathop{\DBcM(\Bs)_F}\nolimits}
\newcommand{\DBcFEM}{\mathop{\DBcM(\Spec E)_F}\nolimits}
\newcommand{\DBcFkM}{\mathop{\DBcM(\Spec k)_F}\nolimits}
\newcommand{\DBcFQM}{\mathop{\DBcM(\BQ)_F}\nolimits}
\newcommand{\DBcQQM}{\mathop{\DBcM(\BQ)_\BQ}\nolimits}

\newcommand{\DBcFUM}{\mathop{\DBcM(U)_F}\nolimits}
\newcommand{\DBcXM}{\mathop{\DBcM(X)}\nolimits}
\newcommand{\DBcFWM}{\mathop{\DBcM(W)_F}\nolimits}
\newcommand{\DBcFXM}{\mathop{\DBcM(X)_F}\nolimits}

\newcommand{\DBcFZM}{\mathop{\DBcM(Z)_F}\nolimits}
\newcommand{\DBcFZoM}{\mathop{\DBcM(W')_F}\nolimits}

\newcommand{\DgM}{\mathop{DM_{gm}(k)}\nolimits}

\newcommand{\DFTBsM}{\mathop{DMT(\Bs)_F}\nolimits}

\newcommand{\ReL}{\mathop{{\rm Res}_{L/\BQ}}\nolimits}
\newcommand{\Mgm}{\mathop{M_{gm}}\nolimits}

\newcommand{\Mcgm}{\mathop{M_{gm}^c}\nolimits}

\newcommand{\dMgm}{\mathop{\partial M_{gm}}\nolimits}
\newcommand{\Bs}{\mathop{B_\sigma}\nolimits}

\newcommand{\Ss}{\mathop{S_\sigma}\nolimits}
\newcommand{\bSs}{\mathop{\overline{\Ss}}\nolimits}

\newcommand{\Zp}{\mathop{Z_\varphi}\nolimits}
\newcommand{\ip}{\mathop{i_\varphi}\nolimits}
\newcommand{\pis}{\mathop{\pi_\sigma}\nolimits}
\newcommand{\pios}{\mathop{\pi'_\sigma}\nolimits}
\newcommand{\pits}{\mathop{\pi''_\sigma}\nolimits}

\newcommand{\one}{\mathds{1}}

\begin{document}

%

\hfuzz=3pt
\overfullrule=10pt                   


\setlength{\abovedisplayskip}{6.0pt plus 3.0pt}
\setlength{\belowdisplayskip}{6.0pt plus 3.0pt}
\setlength{\abovedisplayshortskip}{6.0pt plus 3.0pt}
\setlength{\belowdisplayshortskip}{6.0pt plus 3.0pt}

\setlength{\baselineskip}{13.0pt}
\setlength{\lineskip}{0.0pt}
\setlength{\lineskiplimit}{0.0pt}

%
%

\title{Chow motives without projectivity, II
\forget{
\footnotemark
\footnotetext{To appear in ....}
}
}
\author{\footnotesize by\\ \\
\mbox{\hskip-2cm
\begin{minipage}{6cm} \begin{center} \begin{tabular}{c}
J\"org Wildeshaus \\[0.2cm]
\footnotesize Universit\'e Paris 13\\[-3pt]
\footnotesize Sorbonne Paris Cit\'e \\[-3pt]
\footnotesize LAGA, CNRS (UMR~7539)\\[-3pt]
\footnotesize F-93430 Villetaneuse\\[-3pt]
\footnotesize France\\
{\footnotesize \tt wildesh@math.univ-paris13.fr}
\end{tabular} \end{center} \end{minipage}
\hskip-2cm}
\\[2.5cm]
}
\date{September 27, 2018}
\maketitle
\begin{abstract}
\noindent
The purpose of this article is to provide a simplified construction 
of the intermediate extension of a Chow motive, provided a condition
on absence of weights in the boundary is satisfied. We give a criterion, 
which guarantees the validity of the condition, and 
compare our new construction
to the theory of the interior motive established earlier. 
We finish the article with a review of the known applications 
to the boundary of Shimura varieties. \\

\noindent Keywords: weight structures, 
boundary motive, interior motive, motivic intermediate extension, intersection motive.

\end{abstract}


\bigskip
\bigskip
\bigskip

\noindent {\footnotesize Math.\ Subj.\ Class.\ (2010) numbers: 
14F42
(14C25, 14F20, 14F32, 18E30, 19E15).
}

\eject

\tableofcontents

\bigskip
\vspace*{0.5cm}


%
%

\setcounter{section}{0}
\section{Introduction}
\label{Intro}



The aim of this article is to extend the main results from \cite{W4} to the context
of motives over a base scheme $X$, taking into account and relying on the substantial
progress the motivic theory has undergone since the writing of \loccit. \\

As far as our aim is concerned,
this progress concerns two main points: (1)~the construction of the triangulated category
$\DBcXM$ of motives over $X$ (generalizing the $\BQ$-linear version of Voevodsky's 
definition for $X$ equal to a point, \emph{i.e.}, to the spectrum of a field), together with the \emph{formalism
of six operations}, (2)~the construction of a \emph{weight structure} on $\DBcXM$, compa\-tible
with the six operations. \\   

As in \cite{W4}, the focus of our study is the \emph{absence of weights},
and the guiding principle remains that absence of weights in motives associated to
a \emph{boundary} allows for the construction
of a privileged \emph{extension} of a given (Chow) motive. 
Even over a point, our approach \emph{via} relative motives
yields a new criterion (Theorem~\ref{2N}) on absence of weights in the boundary. \\

In the geometrical context of \emph{Siegel threefolds},
it is that new criterion that is needed to control the weights \cite{W12}.
Indeed, the observation that earlier results concerning motives over a point were
not sufficient to analyze the boundary of Siegel threefolds, can
be seen as the main
motivation of the present paper. \\ 

For a scheme $U$, which is separated and of finite type over
a field $k$ (assumed to be of characteristic zero, to fix ideas), the \emph{boundary
motive} $\dMgm(U)$ of $U$ \cite{W2} fits into an exact triangle 
\[
(\ast) \quad\quad
\dMgm(U) \longto \Mgm(U) \stackrel{u}{\longto} \Mcgm(U) \longto \dMgm(U)[1] \; .
\]
Here, $\Mgm(U)$ and $\Mcgm(U)$ denote the \emph{motive} and the \emph{motive
with compact support}, respectively, of $U$, as defined by Voevodsky \cite{V}. Assuming in addition that $U$
is smooth over $k$, the objects $\Mgm(U)$ and $\Mcgm(U)$ are of weights $\le 0$
and $\ge 0$, respectively. The axioms imposed on a weight structure then formally imply
that the morphism $u$ factors over a motive, which is pure of weight zero,
in other words, $u$ factors over a Chow motive over $k$. 
However, such a factorization is by no means unique (for example, the motive of any smooth 
compactification of $U$ provides such a factorization).
In this context, which is the one studied in \cite{W4},
the ``boundary'' is the boundary motive $\dMgm(U)$,
and any factorization of $u$ through a Chow motive is an ``extension''. \\

The above-mentioned progress, and more particularly, point~(1) allows for what one
might call ``geometrical realizations'' of the exact triangle $(\ast)$. Indeed, any
open immersion $j: U \into X$ gives rise to an exact triangle 
\[
(\ast \ast) \quad\quad
i_*i^*j_* \one_U[-1] \longto j_! \one_U \stackrel{m}{\longto} j_* \one_U \longto 
i_*i^*j_* \one_U
\]
of motives over $X$. Here, we denote by $i$ the closed immersion of the complement $Z$
of $U$ into $X$, and by $\one_U$ the structural motive over $U$. 
Provided that the structure morphism $a$ of $X$ is proper, the direct image $a_*$ of
$(\ast \ast)$, \emph{i.e.}, the exact triangle
\[
a_*i_*i^*j_* \one_U[-1] \longto a_*j_! \one_U \stackrel{a_*m}{\longto} a_*j_* \one_U \longto 
a_*i_*i^*j_* \one_U \; ,
\]
is isomorphic to the dual of $(\ast)$. \\

Thanks to point~(2),
relative motives are endowed with weights. Independently of properness 
of the morphism $a$, the motives $j_! \one_U$ and $j_* \one_U$ over $X$
are of weights $\le 0$ and $\ge 0$, respectively. Again, the axioms for weight structures imply
that $m$ factors over a motive of weight zero.
In this relative context, 
the ``boundary'' is the motive $i_*i^*j_* \one_U$ over $X$ (or equivalently, the motive $i^*j_* \one_U$ over $Z$),
and any factorization of $m$ through a Chow motive is an ``extension'' of $\one_U$ (note that unless $X$ is smooth,
too, the structural motive $\one_X$ is in general not pure of weight zero). 
If $a$ is proper, then the functor $a_*$ is \emph{weight exact}. 
Applying it to a factorization of $m$ through a Chow motive over $X$,
one therefore obtains a factorization of $u$ through a Chow motive over $k$. \\

More generally, any Chow motive
$M_U$ over $U$ yields an exact triangle 
\[
i_*i^*j_* M_U[-1] \longto j_! M_U \stackrel{m}{\longto} j_* M_U \longto 
i_*i^*j_* M_U \; ,
\]
$j_! M_U$ and $j_* M_U$ are of weights $\le 0$ and $\ge 0$, respectively, and therefore, the morphism
$m$ factors over a Chow motive over $X$. 
It is this context of relative motives which seems to be best adapted to our study.
The above mentioned guiding principle relates \emph{absence of weights $0$ and $1$}
in the boundary $i^*j_* M_U$ to the existence of a privileged
extension of $M_U$ to a Chow motive over $X$, which is minimal in a precise sense
among all such extensions. Furthermore, 
we are able to describe the sub-category of Chow motives over $X$
arising as such extensions. Our first main result is Theorem~\ref{0B};
it states that restriction $j^*$ from $X$ to $U$
induces an equivalence of categories
\[
j^* : \CHXM_{i^*w \le -1,i^!w \ge 1} \isoto \CHUM_{\partial w \ne 0,1} \, ,
\]
where the left hand side denotes the full sub-category of Chow motives $M$ over $X$
such that $i^*M$ is of weights at most $-1$, and $i^!M$ is of weights at least $1$,
and the right hand side denotes the full sub-category of Chow motives $M_U$ over $U$
such that $i^*j_* M_U$ is \emph{without weights $0$ and $1$}. \\

It turns out that Theorem~\ref{0B} is best proved in the abstract setting of 
triangulated categories $\CC(U)$, $\CC(X)$ and $\CC(Z)$ related
by \emph{gluing}, and equipped with weight structures \emph{compatible with the gluing}.
This is the setting of Section~\ref{0}. The proof of Theorem~\ref{0B}
relies on Construction~\ref{Cons}, which relates 
factorizations of $m : j_! M_U \to j_* M_U$,
for objects $M_U$ of the \emph{heart} $\CC(U)_{w=0}$ of the weight structure on
$\CC(U)$, to \emph{weight filtrations}
of $i^*j_* M_U$. 
Theorem~\ref{0B} establishes an equivalence
 \[
j^* : \CC(X)_{w=0,i^*w \le -1,i^!w \ge 1} \isoto \CC(U)_{w=0,\partial w \ne 0,1} \, ,
\]
where source and target are defined in obvious analogy with the motivic situation. 
We can thus define the \emph{restriction of the intermediate extension
to the category} $\CC(U)_{w=0,\partial w \ne 0,1}$
\[
\ujast : \CC(U)_{w=0,\partial w \ne 0,1} \longinto \CC(X)_{w=0}
\] 
as the composition of the inverse of the equivalence
of Theorem~\ref{0B}, followed by the inclusion
$\CC(X)_{w=0,i^*w \le -1,i^!w \ge 1} \longinto \CC(X)_{w=0}$ (Definition~\ref{0C}). \\

Theorem~\ref{0B} allows us to provide important complements for the existing theory.
First (Remark~\ref{0D}), the functor $\ujast$ is compatible with
the theory developed in 
\cite[Sect.~2]{W9} when the additional hypothesis enabling the set-up
of the latter, namely \emph{semi-primality} of the category $\CC(Z)_{w=0} \,$, is satisfied.
Note that the functoriality properties of the theory from \loccit \ are intrinsically incomplete
as the target of the intermediate extension is only a quotient of $\CC(X)_{w=0} \, $.
Definition~\ref{0C} can thus be seen as providing a \emph{rigidification} 
of the intermediate extension on the sub-category $\CC(U)_{w=0,\partial w \ne 0,1}$
of $\CC(U)_{w=0} \, $. This observation has rather useful consequences.
When $\CC(Z)_{w=0}$ is semi-primary, then by our very Definition~\ref{0C},
it is possible to read off
$i^* \ujast M_U$ and $i^! \ujast M_U$
whether or not $M_U$ belongs to $\CC(U)_{w=0,\partial w \ne 0,1} \,$:
indeed, $M_U \in \CC(U)_{w=0,\partial w \ne 0,1}$ if and only if
$i^* \ujast M_U$ is of weights at most $-1$, and $i^! \ujast M_U$ of weights at least $1$.
In particular, the non-rigidified intermediate extension $\ujast M_U$
from \cite[Sect.~2]{W9}
is rigid \emph{a posteriori}, if $\ujast M_U$ belongs to the full sub-category $\CC(X)_{w=0,i^*w \le -1,i^!w \ge 1}$
of $\CC(X)_{w=0} \,$.
Furthermore (Theo\-rem~\ref{0E}), for $M_U \in \CC(U)_{w=0,\partial w \ne 0,1} \,$, 
the interval $[\alpha,\beta] \supset [0,1]$ of weights
avoided by $i^*j_* M_U$ can be determined 
directly from $i^* \ujast M_U$ and $i^! \ujast M_U$. 
For example, 
in the context studied in \cite{W12}, \emph{i.e.}, of motives over Siegel threefolds,
the condition on semi-primality is satisfied, and therefore,
Theo\-rem~\ref{0E} applies. \\

Second, in the motivic context, Theorem~\ref{0B} provides
a criterion on absence of weights in the boundary, provided
that the structure morphism $a$ of $X$ is proper. More generally, if $a$ is any proper
morphism with source $X$, and
$M_U \in \CHUM_{\partial w \ne 0,1} \,$, then thanks to weight exactness of $a_*$,
the motive $a_* i_* i^* j_* M_U$ is still without weights $0$ and $1$. This means that
condition \cite[Asp.~2.3]{W4} is satisfied for the morphism
\[
a_* m: a_*j_! M_U \longto a_* j_* M_U \; .
\]
The principal aim of Section~\ref{1} is to spell out the consequences 
for our situation of the general theory
developed in \cite[Sect.~2]{W4}, given the validity of \cite[Asp.~2.3]{W4},
and to relate them to the restriction of the intermediate extension
to $\CHUM_{\partial w \ne 0,1}$ (Theorems~\ref{2G}--\ref{2I}). 
Let us mention Theorem~\ref{2H} in particular: 
any endomorphism of $(a \circ j)_! M_U$ or of $(a \circ j)_* M_U$
induces an endomorphism of the Chow motive $a_* \ujast M_U$. 
Theorem~\ref{2H} applies in particular to endomorphisms 
``of Hecke type''; again,
this general observation is used in particular in the geometrical context
of Siegel threefolds \cite{W12}. In case the proper morphism $a$ equals
the structure morphism of $X$, the Chow motive $a_* \ujast M_U$ is defined to be
the \emph{intersection motive of $U$ relative to $X$ with coefficients in $M_U$}
(Definition~\ref{2Ia}). Given the state of the literature, it appeared
useful to spell out the isomorphism between the dual of the \emph{interior motive}
\cite[Sect.~4]{W4} and the intersection motive. The comparison results from Proposition~\ref{2K}
onwards contain the earlier mentioned isomorphism between the dual of 
the exact triangle 
\[
(\ast) \quad\quad
\dMgm(U) \longto \Mgm(U) \stackrel{u}{\longto} \Mcgm(U) \longto \dMgm(U)[1]
\]
and the exact triangle
\[
a_*i_*i^*j_* \one_U[-1] \longto a_*j_! \one_U \stackrel{a_*m}{\longto} a_*j_* \one_U \longto 
a_*i_*i^*j_* \one_U \; .
\]

At this point, it is probably useful to recall that the quest for a motivic analogue
of \emph{intersection cohomology} started some thirty years ago, with the successful
construction by Scholl of what should nowadays be seen as the intersection motive
of \emph{modular curves} \cite{S}. This example, as other examples concerning 
\emph{Shimura varieties}, will be discussed in Section~\ref{3}. 
To the best of the author's knowledge, the only case where an intersection motive
of ``non-Shimura type''
was constructed over the field of definition of the geometric object, concerns 
arbitrary surfaces (with constant coefficients) \cite{CM}; it may be worthwhile to note that
this result appeared almost fifteen years after Scholl's! \\

A concrete difficulty arises when the defining condition
of $\CHUM_{\partial w \ne 0,1}$ needs to be checked for
a concrete object of $\CHUM$: given a Chow motive $M_U$
over $U$, how to determine whether or not the motive $i^* j_* M_U$
is without weights $0$ and $1$? 
Section~\ref{2} gives what we think of as the optimal answer that can be given to date.
Combining key results from \cite{W11} and \cite{W9}, we prove  
Theorem~\ref{2C}, which we consider as our second main result:
assume that the (generic) \emph{$\ell$-adic realization} $R_{\ell,U}(M_U)$ of $M_U$
is concentrated in a single perverse degree, and that the motive $i^* j_* M_U$ is 
\emph{of Abelian type} \cite{W9}. Then whether or not $M_U$ belongs to $\CHUM_{\partial w \ne 0,1}$
can be read off the perverse cohomology sheaves of $i^* \ujast R_{\ell,U}(M_U)$
and of $i^! \ujast R_{\ell,U}(M_U)$. If $M_U \in \CHUM_{\partial w \ne 0,1}$,
then the precise interval 
$[\alpha,\beta] \supset [0,1]$ of weights
avoided by $i^*j_* M_U$ can be determined 
from $i^* \ujast R_{\ell,U}(M_U)$ and $i^! \ujast R_{\ell,U}(M_U)$. \\

Chow motives have a tendancy to be auto-dual up to a shift and a twist; this
is in any case true for the Chow motives occurring in the applications we have in mind,
\emph{e.g.}, in the earlier mentioned analysis of the weights 
in the boundary of Siegel threefolds \cite{W12}. 
Given that the criterion from Theorem~\ref{2C} is compatible with duality,
one may hope that the verification of a certain half of that criterion,
for example the half concerning $i^* \ujast R_{\ell,U}(M_U)$,
is sufficient, when $M_U$ is auto-dual. This hope is made explicit in   
Corollary~\ref{2Ca}.
We think of this result as potentially very useful for other applications. 
For the sake of completeness, we combine Corollary~\ref{2Ca} with the comparison 
from Section~\ref{1}; the result is the earlier mentioned Theorem~\ref{2N}.  \\

The final Section~\ref{3} contains a complete review of the known applications of our
theory to Shimura varieties. Let us point out that some of these cases are equally covered by 
a recent result of Vaish's \cite{Va}. His approach replaces weight structures
by weight truncations \`a la S.~Morel (but still relies on the main result from \cite{W11}),
thereby providing an alternative approach to the problem of rigidification 
of the intermediate extension. It is interesting to note that Vaish's result applies in 
certain situations where our condition on absence of weights $0$ and $1$
is not satisfied. \\

Part of this work was done while I was enjoying a \emph{d\'el\'egation
aupr\`es du CNRS}, to which I wish to express my gratitude. 
I also wish to thank F.~D\'eglise for useful discussions,
and the referees for their comments, which contribued 
considerably to improve this article.\\

{\bf Conventions}: Throughout the article, 
$F$ denotes a finite direct product of fields
of characteristic zero, in other words, a commutative semi-simple Noetherian $\BQ$-algebra.
We fix a base scheme
$\BB$, which is of finite type over some excellent scheme 
of dimension at most two. By definition, \emph{schemes} are  
$\BB$-schemes which are separated and 
of finite type (in particular, they are excellent, and
Noetherian of finite dimension), \emph{morphisms} between schemes
are separated morphisms of $\BB$-schemes, and 
a scheme is \emph{nilregular}
if the underlying reduced scheme is regular. \\

We use the triangulated, $\BQ$-linear categories
$\DBcXM$ of \emph{constructible Beilinson motives} over $X$ 
\cite[Def.~15.1.1]{CD},  
indexed by schemes $X$ (always in the sense of the above conventions). 
In order to have an $F$-linear theory at one's disposal, 
one re-does the construction, but using $F$
instead of $\BQ$ as coefficients \cite[Sect.~15.2.5]{CD}. This yields 
triangulated, $F$-linear categories
$\DBcFXM$ satisfying the $F$-linear analogues of the properties
of $\DBcXM$. In particular, these categories 
are pseudo-Abelian (see \cite[Sect.~2.10]{H}).
Furthermore, the canonical functor of $F$-linear categories $\DBcXM \otimes_\BQ F \to \DBcFXM$
is fully faithful \cite[Sect.~14.2.20]{CD}.
As in \cite{CD}, the symbol $\one_X$
is used to denote the unit for the tensor product in $\DBcFXM$.
We shall employ the full formalism of six operations developed in
\loccit . The reader may choose to consult \cite[Sect.~2]{H} or
\cite[Sect.~1]{W10} for concise presentations of this formalism. \\

Beilinson motives can be endowed with a canonical weight structure,
thanks to the main results from \cite{H}
(see \cite[Prop.~6.5.3]{Bo} for the case
$X = \Spec k$, for a field $k$ of characteristic zero). We refer to it as the \emph{motivic
weight structure}. Following \cite[Def.~1.5]{W10}, the category 
$\CHFXM$ of \emph{Chow motives} over $X$
is defined as the heart $\DBcM(X)_{F,w = 0}$ of the motivic weight structure
on $\DBcFXM$. It equals the pseudo-Abelian completion
of $\CHXM_\BQ \otimes_\BQ F$. According to \cite[Thm.~3.3~(ii)]{H},
the motivic weight structure on $\DBcFXM$ is uniquely determined by the requirement that 
$f_*\one_Y (n)[2n] \in \CHFXM$ whenever $n \in \BZ$, and $f: Y \to X$ is a proper morphism
with regular source $Y$. \\

When we assume a field $k$ to \emph{admit resolution of singularities},
then it will be 
in the sense of \cite[Def.~3.4]{FV}:
(i)~for any separated $k$-scheme $X$ of finite type, there exists an abstract blow-up $Y \to X$
\cite[Def.~3.1]{FV} whose source $Y$ is smooth over $k$,
(ii)~for any pair of smooth, seperated $k$-schemes $X, Y$ of finite type, 
and any abstract blow-up $q : Y \to X$,
there exists a sequence of blow-ups 
$p: X_n \to \ldots \to X_1 = X$ with smooth centers,
such that $p$ factors through $q$. 
We say that $k$ \emph{admits strict resolution of singularities},
if in (i), for any given dense open subset $U$
of the smooth locus of $X$,
the blow-up $q: Y \to X$ can be chosen to be an isomorphism above $U$,
and such that arbitrary intersections of
the irreducible components of the complement $Z$ of $U$ in $Y$  
are smooth (e.g., $Z \subset Y$ a normal
crossing divisor with smooth irreducible components).


\bigskip
%
%

\section{Rigidification of the intermediate extension}
\label{0}



Throughout this section, let us fix
three $F$-linear pseudo-Abelian triangulated 
categories $\CC(U)$, $\CC(X)$ and $\CC(Z)$,
the second of which is obtained from the others
by gluing. This means that the three categories are
equipped with six exact functors
\[
\vcenter{\xymatrix@R-10pt{
        \CC(U) \ar[r]^-{j_!} & \CC(X) \ar[r]^-{i^*} & \CC(Z) \\
        \CC(U) & \CC(X) \ar[l]_-{j^*} & \CC(Z) \ar[l]_-{i_*} \\ 
        \CC(U) \ar[r]^-{j_*} & \CC(X) \ar[r]^-{i^!} & \CC(Z)
\\}}
\] 
satisfying the axioms from \cite[Sect.~1.4.3]{BBD}. We assume that
$\CC(U)$, $\CC(X)$ and $\CC(Z)$ are equipped with weight structures
$w$ (the same letter for the three weight structures), and that
the one on $\CC(X)$ is actually obtained from the two others in a
way compatible with the gluing, meaning 
that the left adjoints $j_!$, $j^*$, $i^*$
and $i_*$ respect the categories $\CC(\bullet)_{w \le 0}$, and the right
adjoints $j^*$, $j_*$, $i_*$ and $i^!$ respect the categories
$\CC(\bullet)_{w \ge 0}$. In particular, we have a fully faithful functor
\[
i_* : \CC(Z)_{w = 0} \longinto \CC(X)_{w = 0} \; ,
\]
and a functor 
\[
j^* : \CC(X)_{w = 0} \longto \CC(U)_{w = 0} \; .
\]
According to \cite[Prop.~2.5]{W9}, the latter is full and essentially
surjective. We shall need to understand its restriction to a certain
sub-category of $\CC(X)_{w = 0}$. \\

Recall \cite[Def.~1.10]{W4} that an object $M$ is said to be  
\emph{without weights} $m,\ldots,n$, 
or \emph{to avoid weights} $m,\ldots,n$, for integers $m \le n$, if it admits a 
\emph{weight filtration avoiding weights} $m,\ldots,n$, \emph{i.e.}
\cite[Def.~1.6]{W4}, if there is an exact triangle
\[ 
M_{\le {m-1}} \longto M \longto M_{\ge {n+1}} \longto M_{\le m-1}[1]
\]
with $M_{\le {m-1}}$ of weights at most $m-1$,
and $M_{\ge {n+1}}$ of weights at least $n+1$. 

\begin{Def} \label{0A}
(a)~Denote by $\CC(U)_{w=0,\partial w \ne 0,1}$ 
the full sub-category of $\CC(U)_{w=0}$
of objects $M_U$ such that $i^*j_* M_U$ is without weights $0$ and $1$. 
\\[0.1cm]
(b)~Denote by $\CC(X)_{w=0,i^*w \le -1,i^!w \ge 1}$ 
the full sub-category of $\CC(X)_{w=0}$
of objects $M$ such that $i^* M \in \CC(Z)_{w \le -1}$ and 
$i^! M \in \CC(Z)_{w \ge 1}$.
\end{Def}

\begin{Thm} \label{0B}
The restriction of $j^*$ to $\CC(X)_{w=0,i^*w \le -1,i^!w \ge 1}$
induces an equivalence of categories
\[
j^* : \CC(X)_{w=0,i^*w \le -1,i^!w \ge 1} \isoto 
\CC(U)_{w=0,\partial w \ne 0,1} \; .
\]
\end{Thm}

The proof of Theorem~\ref{0B} relies on the following.

\begin{Cons} \label{Cons}
Let $M_U \in \CC(U)_{w=0}$, and consider the morphism
\[
m: j_! M_U \longto j_* M_U
\]
in $\CC(X)$.
There is a canonical choice of cone of $m$, namely, the object
$i_* i^* j_* M_U$.
Any weight filtration of $i^* j_* M_U$
\[
C_{\le 0} \stackrel{c_-}{\longto} i^* j_* M_U 
\stackrel{c_+}{\longto} C_{\ge 1} 
\stackrel{\delta}{\longto} C_{\le 0}[1]
\]
(with $C_{\le 0} \in \CC(Z)_{w \le 0}$ and $C_{\ge 1} \in \CC(Z)_{w \ge 1}$)
yields a diagram, 
whose $4$-term lines and columns are exact triangles, and
which we shall denote by the symbol~(1)
\[
\vcenter{\xymatrix@R-10pt{
        0 \ar[d] \ar[r] &
        i_* C_{\ge 1}[-1] \ar@{=}[r] &
        i_* C_{\ge 1}[-1] \ar[d]^{i_* \delta[-1]} \ar[r] &
        0 \ar[d] \\
        j_! M_U \ar@{=}[d] &  
         &
        i_* C_{\le 0} \ar[d]^{i_* c_-} \ar[r] &
        j_! M_U[1] \ar@{=}[d] \\
        j_! M_U \ar[d] \ar[r]^m &
        j_* M_U \ar[d] \ar[r] &
        i_* i^* j_* M_U \ar[d]^{i_* c_+} \ar[r] &
        j_! M_U[1] \ar[d] \\
        0 \ar[r] &
        i_* C_{\ge 1} \ar@{=}[r] &
        i_* C_{\ge 1} \ar[r] &
        0    
\\}}
\]
According to axiom TR4' of triangulated categories 
\cite[Sect.~1.1.6]{BBD}, diagram~(1)
can be completed to give a diagram, denoted by~(2)
\[
\vcenter{\xymatrix@R-10pt{
        0 \ar[d] \ar[r] &
        i_* C_{\ge 1}[-1] \ar[d] \ar@{=}[r] &
        i_* C_{\ge 1}[-1] \ar[d]^{i_* \delta[-1]} \ar[r] &
        0 \ar[d] \\
        j_! M_U \ar@{=}[d] \ar[r] &  
        M \ar[d] \ar[r] &
        i_* C_{\le 0} \ar[d]^{i_* c_-} \ar[r] &
        j_! M_U[1] \ar@{=}[d] \\
        j_! M_U \ar[d] \ar[r]^m &
        j_* M_U \ar[d] \ar[r] &
        i_* i^* j_* M_U \ar[d]^{i_* c_+} \ar[r] &
        j_! M_U[1] \ar[d] \\
        0 \ar[r] &
        i_* C_{\ge 1} \ar@{=}[r] &
        i_* C_{\ge 1} \ar[r] &
        0    
\\}}
\] 
with $M \in \CC(X)$. By the second row, and the second column of
diagram~(2), the object $M$
is simultaneously an extension of objects of weights $\le 0$, 
and an extension of objects of weights $\ge 0$.
It follows easily (\emph{c.f.}\ \cite[Prop.~1.3.3~3]{Bo})
that $M$ belongs to both $\CC(X)_{w \le 0}$ and $\CC(X)_{w \ge 0} \, $, and hence to
$\CC(X)_{w=0} \, $.

Applying the functors $j^*$, $i^*$, and $i^!$ to~(2), we see that $j^* M = M_U$, 
\[
i^* M \isoto C_{\le 0} \; , \quad \text{and} \quad
C_{\ge 1}[-1] \isoto i^! M \; .
\]
\end{Cons}

\medskip

\begin{Proofof}{Theorem~\ref{0B}}
For $M_U \in \CC(U)_{w=0}$
and a weight filtration 
\[
C_{\le 0} \longto i^* j_* M_U \longto C_{\ge 1} \longto C_{\le 0}[1]
\] 
of $i^* j_* M_U$, apply Construction~\ref{Cons}
to get an extension $M \in \CC(X)_{w=0}$ of $M_U$ to $X$.
From the isomorphisms
\[
i^* M \isoto C_{\le 0} \quad \text{and} \quad
C_{\ge 1}[-1] \isoto i^! M \; ,
\]
we see that $M \in \CC(X)_{w=0,i^*w \le -1,i^!w \ge 1}$ if and only if
\[
C_{\le 0} \in \CC(Z)_{w \le -1} \quad \text{and} \quad
C_{\ge 1} \in \CC(Z)_{w \ge 2} \; .
\]
In particular, we see that any object in 
$\CC(U)_{w=0,\partial w \ne 0,1}$ admits a pre-image under $j^*$ in
$\CC(X)_{w=0,i^*w \le -1,i^!w \ge 1} \,$.

Conversely, any object $M$ from $\CC(X)_{w=0}$
fits into a diagram of type~(2)
\[
\vcenter{\xymatrix@R-10pt{
        0 \ar[d] \ar[r] &
        i_* i^! M \ar[d] \ar@{=}[r] &
        i_* i^! M \ar[d] \ar[r] &
        0 \ar[d] \\
        j_! j^* M \ar@{=}[d] \ar[r] &  
        M \ar[d] \ar[r] &
        i_* i^* M \ar[d] \ar[r] &
        j_! j^* M[1] \ar@{=}[d] \\
        j_! j^* M \ar[d] \ar[r] &
        j_* j^* M \ar[d] \ar[r] &
        i_* i^* j_* j^* M \ar[d] \ar[r] &
        j_! j^* M[1] \ar[d] \\
        0 \ar[r] &
        i_* i^! M[1] \ar@{=}[r] &
        i_* i^! M[1] \ar[r] &
        0    
\\}}
\]  
Its third column shows that
$j^*$ maps $\CC(X)_{w=0,i^*w \le -1,i^!w \ge 1}$ to 
$\CC(U)_{w=0,\partial w \ne 0,1}$. 

The restriction of $j^*$ therefore yields a well-defined functor
\[
j^* : \CC(X)_{w=0,i^*w \le -1,i^!w \ge 1} \longto 
\CC(U)_{w=0,\partial w \ne 0,1} \; ,
\]
which is full \cite[Prop.~2.5]{W9}
and essentially surjective. It remains to show that it
is faithful, \emph{i.e.}, that a morphism $f$ in 
$\CC(X)_{w=0,i^*w \le -1,i^!w \ge 1}$ such that $j^* f = 0$, equals zero.

Thus, let $M, M' \in \CC(X)_{w=0,i^*w \le -1,i^!w \ge 1}$, $f: M \to M'$,
and assume that $j^* f = 0$. This means that the composition
of $f$ with the adjunction morphism $M' \to j_*j^* M'$ is zero. Given
the exact localization triangle
\[
i_*i^! M' \longto M' \longto j_*j^* M' \longto i_*i^! M' [1] \; ,
\]
the morphism $f$ factors through $i_*i^! M'$. Now $M$ is of weight zero,
while $i^! M'$, and hence $i_*i^! M'$, is of strictly positive weights.
By orthogonality, any morphism $M \to i_*i^! M'$ is zero.
\end{Proofof}

\begin{Def} \label{0C}
Let the $F$-linear pseudo-Abelian
triangulated catego\-ries $\CC(U)$, $\CC(X)$ and $\CC(Z)$
be related by gluing, and equipped with weight structures $w$
compatible with the gluing. 
Define the \emph{restriction of the intermediate extension
to the category} $\CC(U)_{w=0,\partial w \ne 0,1}$
\[
\ujast : \CC(U)_{w=0,\partial w \ne 0,1} \longinto \CC(X)_{w=0}
\] 
as the composition of the inverse of the equivalence
of Theorem~\ref{0B}, followed by the inclusion
$\CC(X)_{w=0,i^*w \le -1,i^!w \ge 1} \longinto \CC(X)_{w=0}$. 
\end{Def}

\begin{Rem} \label{0Ca}
Assume that contravariant auto-equivalences
\[
\BD_\argdot : \CC(\argdot) \isoto \CC(\argdot) \; ,
\argdot \in \{ U,X, Z \}
\]
are given, that they are compatible with the gluing 
(\emph{e.g.}, $\BD_X \circ j_* \cong j_! \circ \BD_U$)
and with the weight structures (\emph{e.g.}, 
$\BD_X \CC(X)_{w \ge 0} = \CC(X)_{w \le 0})$. 
From the isomorphisms $\BD_Z \circ i^* \circ j_* \cong i^! \circ j_! \circ \BD_U$ and
$i^! \circ j_! \cong i^* \circ j_* [-1]$, it follows first that the functor
$\BD_U$ respects the category $\CC(U)_{w=0,\partial w \ne 0,1} \,$.
Similarly, the functor
$\BD_X$ respects the category $\CC(X)_{w=0,i^*w \le -1,i^!w \ge 1} \,$.

It then follows formally from Definition~\ref{0C},
and from $\BD_U \circ j^* \cong j^* \circ \BD_X$ that 
\[
\BD_X \circ j_{!*} \cong j_{!*} \circ \BD_U \; .
\]
In other words, the restriction of the intermediate extension
to $\CC(U)_{w=0,\partial w \ne 0,1}$
is compatible with local duality.
\end{Rem}

\begin{Rem} \label{0D}
(a)~Assume that some composition of morphisms
\[
i_* N \longto M \longto i_* N 
\]
gives the identity on $i_* N$,
for some object $N$ of $\CC(Z)_{w = 0}$. Then
the adjunction properties of $i^*$, $i_*$ and $i^!$,
and in particular, the identifications $i^*i_* = \id_{\CC(Z)}$
and $\id_{\CC(Z)} = i^! i_*$,
show that $N$ is a direct factor of both
$i^! M$ and $i^* M$, 
and that the restriction of the composition $i^! M \to i^* M$
of the adjunction morphisms to this direct factor is the identity. 

We obtain that objects in $\CC(X)_{w=0,i^*w \le -1,i^!w \ge 1}$
do not admit non-zero direct factors belonging to the image of
$\CC(Z)_{w=0}$ under the functor $i_*$. This justifies Definition~\ref{0C}:
the intermediate extension of $M_U \in \CC(U)_{w=0,\partial w \ne 0,1}$
is indeed an extension of $M_U$ not admitting non-zero direct factors belonging 
to the image of $i_*$. \\[0.1cm]
(b)~More generally, for $M \in \CC(X)_{w=0,i^*w \le -1,i^!w \ge 1}$
and $N \in \CC(Z)_{w=0} \, $, any morphism $M \to i_* N$ is zero, and so is any
morphism $i_* N \to M$. \\[0.1cm]
(c)~Following \cite[Def.~1.6~(a)]{W9},
denote by $\Fg$ the two-sided ideal of 
$\CC(X)_{w=0}$ generated by
\[
\Hom_{\CC(X)_{w=0}} (M,i_*N) \quad \text{and} \quad 
\Hom_{\CC(X)_{w=0}} (i_*N,M) \; ,
\]
for all objects $(M,N)$ of 
$\CC(X)_{w=0} \times 
\CC(Z)_{w=0} \, $, 
such that $M$ admits no non-zero direct factor belonging to the image of
$i_*$. Denote by $\CC(X)_{w=0}^u$ 
the quotient of $\CC(X)_{w=0}$ by $\Fg$.
From (b), we see that $\Fg(M,M') = 0$, for any 
$M, M' \in \CC(X)_{w=0,i^*w \le -1,i^!w \ge 1}$. Thus, the 
quotient functor $\CC(X)_{w=0} \onto \CC(X)_{w=0}^u$
induces an auto-equivalence on $\CC(X)_{w=0,i^*w \le -1,i^!w \ge 1}$. \\[0.1cm]
(d)~From the above, we see that Definition~\ref{0C} is compatible with
the theory developed in 
\cite[Sect.~2]{W9} when the hypotheses enabling the set-up
of the latter are satisfied.
More precisely, assume in addition that $\CC(Z)_{w=0}$ is semi-primary
\cite[D\'ef.~2.3.1]{AK}. Then 
\[
\ujast : \CC(U)_{w=0} \longinto \CC(X)_{w=0}^u
\] 
is defined on the whole of $\CC(U)_{w=0}$ \cite[Def.~2.10]{W9}.
Parts~(a) and (c) of the present Remark, and \cite[Summ.~2.12~(a)~(1)]{W9} show
that the diagram
\[
\vcenter{\xymatrix@R-10pt{
        \CC(U)_{w=0,\partial w \ne 0,1} 
               \ar@{_{ (}->}[d] \ar@{^{ (}->}[r]^-{\ujast} &
        \CC(X)_{w=0}
               \ar@{->>}[d] \\
        \CC(U)_{w=0} \ar@{^{ (}->}[r]^-{\ujast} &
        \CC(X)_{w=0}^u         
\\}}
\] 
is commutative.
Thus, $\ujast : \CC(U)_{w=0,\partial w \ne 0,1} \into \CC(X)_{w=0}$ 
from Definition~\ref{0C}
is indeed the restriction of $\ujast : \CC(U)_{w=0} \into \CC(X)_{w=0}^u$
from \cite[Sect.~2]{W9},
when the latter is defined, to $\CC(U)_{w=0,\partial w \ne 0,1}$.
\end{Rem}

\begin{Thm} \label{0E}
Let $M_U \in \CC(U)_{w=0} \, $. \\[0.1cm]
(a)~Assume that the category $\CC(Z)_{w=0}$ is semi-primary,
so that the functor $\ujast$ is defined on the whole of
$\CC(U)_{w=0}$ \cite[Sect.~2]{W9}. 
Then the object $M_U$ belongs to $\CC(U)_{w=0,\partial w \ne 0,1}$ if and only if
$\ujast M_U$ belongs to $\CC(X)_{w=0,i^*w \le -1,i^!w \ge 1} \, $. \\[0.1cm]
(b)~Assume that $M_U$ belongs to $\CC(U)_{w=0,\partial w \ne 0,1}$.
Let $\alpha \le 0$ and $\beta \ge 1$ be two integers. Then $i^* j_* M_U$ is without weights
$\alpha, \alpha+1, \ldots, \beta$ if and only if
\[
i^* \ujast M_U \in \CC(Z)_{w \le \alpha - 1} \quad \text{and} \quad
i^! \ujast M_U \in \CC(Z)_{w \ge \beta} \; .
\]
\end{Thm}

\begin{Proof}
The ``if'' part from (a) is Theorem~\ref{0B}, and its ``only if'' part is Remark~\ref{0D}~(d). 

As for (b), note that by definition,
Construction~\ref{Cons} for $M_U$ and a weight filtration 
\[
C_{\le -1} \longto i^* j_* M_U \longto C_{\ge 2} \longto C_{\le -1}[1]
\] 
of $i^* j_* M_U$
avoiding weights $0$ and $1$ yields the extension $M = \ujast M_U$
of $M_U$ to $X$. The claim thus follows from the isomorphisms
\[
i^* \ujast M_U \isoto C_{\le -1} \quad \text{and} \quad
C_{\ge 2}[-1] \isoto i^! \ujast M_U \; .
\]
\end{Proof}

\begin{Rem} \label{0F}
(a)~As the attentive reader will have remarked already, the formalism could have been
developed on larger sub-categories of $\CC(X)$ and $\CC(U)$, at the cost of losing
its inherent auto-duality. More precisely, define full sub-categories  
$\CC(U)_{w=0,\partial w \ne 0}$ and $\CC(U)_{w=0,\partial w \ne 1}$ 
of $\CC(U)_{w=0} \,$, and $\CC(X)_{w=0,i^*w \le -1}$ and $\CC(X)_{w=0,i^!w \ge 1}$ 
of $\CC(X)_{w=0}$ in the obvious way. Then as in Theorem~\ref{0B},
\[
j^* : \CC(X)_{w=0,i^*w \le -1} \isoto 
\CC(U)_{w=0,\partial w \ne 0} 
\]
and
\[
j^* : \CC(X)_{w=0,i^!w \ge 1} \isoto 
\CC(U)_{w=0,\partial w \ne 1} \; ,
\]
allowing to define the restrictions $\ujast$ of the intermediate extension to both
$\CC(U)_{w=0,\partial w \ne 0}$ and $\CC(U)_{w=0,\partial w \ne 1}$.
Local duality as in Remark~\ref{0Ca} exchanges the two constructions.
Parts~(a), (c), and (d) (but not (b)) of Remark~\ref{0D} apply
with the most obvious modifications. \\[0.1cm]
(b)~Even if  it will not be needed in the sequel of this article,
it is worthwhile to spell out the modified version of Theorem~\ref{0E}.
Let $M_U \in \CC(U)_{w=0} \, $. Then  
\[
M_U \in \CC(U)_{w=0,\partial w \ne 0} \Longleftrightarrow 
\ujast M_U \in \CC(X)_{w=0,i^*w \le -1} 
\]
and 
\[
M_U \in \CC(U)_{w=0,\partial w \ne 1} \Longleftrightarrow 
\ujast M_U \in \CC(X)_{w=0,i^!w \ge 1} \; ,
\]
provided that $\CC(Z)_{w=0}$ is semi-primary.
More interestingly, part~(b) of Theorem~\ref{0E} can be separated into
two statements: let $\alpha \le 0$ and $\beta \ge 1$. 
Assume that $M_U \in \CC(U)_{w=0,\partial w \ne 0}$ or $M_U \in \CC(U)_{w=0,\partial w \ne 1}$.
Then $i^* j_* M_U$ is without weights
$\alpha, \alpha+1, \ldots, 0$ if and only if
\[
i^* \ujast M_U \in \CC(Z)_{w \le \alpha - 1} \; ,
\]
and $i^* j_* M_U$ is without weights
$1, 2, \ldots, \beta$ if and only if
\[
i^! \ujast M_U \in \CC(Z)_{w \ge \beta} \; .
\]
\end{Rem}


\bigskip
%
%

\section{Motivic intermediate extension and interior motive}
\label{1}



The purpose of the present section is to 
connect Section~\ref{0} to the theory developed in \cite[Sect.~2]{W4}. 
The main results are Theorems~\ref{2G}, \ref{2H} and \ref{2I}; these 
concern the \emph{motivic intermediate extension}, and are formal
analogues of the main results from \cite[Sect.~4]{W4} on the \emph{interior motive},
defined and studied in \loccit. When the base scheme $\BB$ is the spectrum of a field
admitting strict resolution of singularities, then the analogy is not just formal:
Corollary~\ref{2L} establishes an isomorphism between the dual of the interior motive and
the direct image of the motivic intermediate extension under the structure morphism,
provided the latter is proper and that weights $0$ and $1$ are avoided. \\
 
Let $X$ be a scheme (in the sense of our Introduction), and $j: U \into X$
an open immersion with complement $i: Z \into X$. 
Thanks to \emph{localization} \cite[Sect.~2.3]{CD},
and to compatibility of the motivic weight structure with gluing \cite[Thm.~3.8]{H},
Theorem~\ref{0B} applies to the categories $\CC(\bullet) = \DBcFbM$, for $\bullet = U, X, Z$.
We write $\CHUM_{F,\partial w \ne 0,1}$ for $\DBcM(U)_{F,w=0,\partial w \ne 0,1} \, $,
and $\CHXM_{F,i^*w \le -1,i^!w \ge 1}$ for $\DBcM(X)_{F,w=0,i^*w \le -1,i^!w \ge 1}$.
In this motivic context, Remark~\ref{0D}~(d) says that the functor $\ujast$ 
from Definition~\ref{0C} equals the restriction 
to $\CHUM_{F,\partial w \ne 0,1}$ of the motivic intermediate extension,
when the latter is defined; note that this is the case 
in the context of \emph{motives of Abelian type}, 
studied in \cite[Sect.~5]{W9} (see our Section~\ref{2}).

\begin{Prop} \label{2D}
Let $a: X \to \BA$ a proper morphism of schemes, and
$M_U$ a Chow motive over $U$ belonging to $\CHUM_{F,\partial w \ne 0,1} \, $.
Denote by $m$ the morphism $j_! M_U \to j_* M_U$. \\[0.1cm]
(a)~The image 
\[
a_* j_! M_U \stackrel{a_* m}{\longto} a_* j_* M_U \longto a_* i_* i^* j_* M_U
\longto a_* j_! M_U [1]
\]
under $a_*$ of the triangle
\[
j_! M_U \stackrel{m}{\longto} j_* M_U \longto i_* i^* j_* M_U
\longto j_! M_U [1]
\]
in $\DBcFDM$ is exact. \\[0.1cm]
(b)~The motive $a_* i_* i^* j_* M_U \in \DBcFDM$ is without weights $0$ and $1$.
\end{Prop}

\begin{Proof}
(a): Indeed, the triangle 
\[
j_! M_U \stackrel{m}{\longto} j_* M_U \longto i_* i^* j_* M_U
\longto j_! M_U [1]
\]
is exact.

(b): The morphism $a \circ i$ being proper, we have 
$(a \circ i)_! = (a \circ i)_* = a_*i_*$. Thus,
$a_*i_*$ is weight exact. It therefore
transforms any weight filtration avoiding weights $0$ and $1$ 
into the same kind of weight filtration.
\end{Proof} 

Thus, \cite[Asp.~2.3]{W4} is satisfied, with $\CC = \DBcFDM$, 
$u = a_* m: a_* j_! M_U \to a_* j_* M_U$, and $C = a_* i_* i^* j_* M_U [-1]$.
Consequently, the theory developed in \cite[Sect.~2]{W4} applies.

\begin{Def}[{cmp.~\cite[Def.~2.1]{W4}}] \label{2E}
Denote by $\DBcDM_{F,w \le 0,\ne -1}$ 
the full sub-category of $\DBcDM_{F,w \le 0}$ of objects without weight $-1$, 
and by $\DBcDM_{F,w \ge 0,\ne 1}$ 
the full sub-category of $\DBcDM_{F, w \ge 0}$ of objects without weight $1$.
\end{Def}

\begin{Prop}[{\cite[Prop.~2.2]{W4}}] \label{2F}
The inclusions 
\[
\iota_-: \CHFDM \longinto \DBcDM_{F, w \le 0,\ne -1} 
\]
and
\[ 
\iota_+: \CHFDM \into \DBcDM_{F, w \ge 0,\ne 1}
\] 
admit a left adjoint 
\[
\Gr_0: \DBcDM_{F, w \le 0,\ne -1} \longto \CHFDM  
\]
and a right adjoint
\[
\Gr_0: \DBcDM_{F, w \ge 0,\ne 1} \longto \CHFDM \; ,
\]
respectively.
Both adjoints map objects (and morphisms) to the term of weight zero 
of a weight filtration avoiding weight $-1$ and $1$, respectively.
The compositions $\Gr_0 \circ \iota_-$ and $\Gr_0 \circ \iota_+$
both equal the identity on $\CHFDM$.
\end{Prop}

\begin{Thm} \label{2G}
Let $a: X \to \BA$ a proper morphism. \\[0.1cm]
(a)~The essential image of the restriction of the 
functor $a_* j_!$ to the sub-category 
$\CHUM_{F,\partial w \ne 0,1}$
is contained in 
$\DBcDM_{F, w \le 0,\ne -1}$, inducing a functor 
\[
(a \circ j)_! = a_* j_! : 
\CHUM_{F,\partial w \ne 0,1}
\longto \DBcDM_{F, w \le 0,\ne -1} \; .
\]
More precisely, if $M_U \in \CHUM_F$ 
is such that
$i^* j_* M_U$ avoids weights $\alpha, \alpha + 1, \ldots, \beta$, for integers
$\alpha \le 0$ and $\beta \ge 1$, then 
\[
a_*i_*i^* \ujast M_U[-1] \longto a_*j_! M_U  \longto a_* \ujast M_U \longto
a_*i_*i^* \ujast M_U
\]
is a weight filtration of $a_*j_! M_U$ avoiding weights
$\alpha-1, \alpha, \ldots, -1$. \\[0.1cm]
(b)~The essential image of the restriction of the 
functor $a_* j_*$ to the sub-category 
$\CHUM_{F,\partial w \ne 0,1}$
is contained in 
$\DBcDM_{F, w \ge 0,\ne 1}$, inducing a functor 
\[
(a \circ j)_* = a_* j_* : 
\CHUM_{F,\partial w \ne 0,1}
\longto \DBcDM_{F, w \ge 0,\ne 1} \; .
\]
More precisely, if $M_U \in \CHUM_F$ 
is such that
$i^* j_* M_U$ avoids weights $\alpha, \alpha + 1, \ldots, \beta$, for integers
$\alpha \le 0$ and $\beta \ge 1$, then 
\[
a_* \ujast M_U \longto a_* j_*M_U \longto a_* i_*i^! \ujast M_U[1] \longto 
a_* \ujast M_U[1] 
\]
is a weight filtration of $a_*j_* M_U$ avoiding weights
$1, 2, \ldots, \beta$. \\[0.1cm]
(c)~There are canonical isomorphisms of functors 
\[
\Gr_0 \circ (a \circ j)_! = \Gr_0 \circ a_* j_! \isoto a_* \ujast 
\quad \text{and} \quad 
a_* \ujast \isoto \Gr_0 \circ a_* j_* = \Gr_0 \circ (a \circ j)_*
\]
on the category
$\CHUM_{F,\partial w \ne 0,1}$. Their composition equals the value  
of the functor
$\Gr_0 \circ a_*$ at the restriction of the natural transformation $m : j_! \to j_*$ to
$\CHUM_{F,\partial w \ne 0,1}$; in particular, 
\[
\Gr_0 \circ (a_* \circ m): \Gr_0 \circ (a \circ j)_! \longto \Gr_0 \circ (a \circ j)_*
\] 
is an isomorphism of functors on $\CHUM_{F,\partial w \ne 0,1}$.  
\end{Thm}

\begin{Proof}
Let $M_U \in \CHUM_{F,\partial w \ne 0,1} \, $. 
By definition (and Theorem~\ref{0B}), the motive $\ujast M_U$ belongs to 
$\CHXM_{F,i^*w \le -1,i^!w \ge 1}$.
Thus, the exact triangles
\[
i_*i^* \ujast M_U[-1] \longto j_! M_U  \longto \ujast M_U \longto
i_*i^* \ujast M_U
\]
and
\[
\ujast M_U \longto j_*M_U \longto i_*i^! \ujast M_U[1] \longto 
\ujast M_U[1] 
\]
are weight filtrations (of $j_! M_U$) avoiding weight $-1$,
and (of $j_*M_U$) avoiding weight $1$, respectively. 
An analogous statement is therefore true for their images under 
the weight exact functor $a_*$ (recall that $a$ is assumed to be proper).
Together with Proposition~\ref{2F},
this shows part~(c) of the statement, as well as the first claims of
parts~(a) and (b). The second, more precise claims follow from Theorem~\ref{0E}~(b).
\end{Proof}

At first sight, it may thus appear that the theory from \cite[Sect.~2]{W4}
does not add much to what we get by explicit identification of the
weight filtrations. But then, note the following.

\begin{Thm} \label{2H}
Let $a: X \to \BA$ a proper morphism. Then for any Chow motive
$M_U \in \CHUM_{F,\partial w \ne 0,1} \, $,
the Chow motive $a_* \ujast M_U \in \CHFDM$ behaves
functorially with respect to both motives $(a \circ j)_! M_U$ and 
$(a \circ j)_* M_U$.
In particular, any endomorphism of $(a \circ j)_! M_U$ or of $(a \circ j)_* M_U$
induces an endomorphism of $a_* \ujast M_U$. 
\end{Thm}

\begin{Proof}
This follows from the functorial identities from Theorem~\ref{2G}~(c).
\end{Proof}

\begin{Thm} \label{2I}
Let $M_U \in \CHUM_{F,\partial w \ne 0,1} \, $, 
$a: X \to \BA$ a proper morphism, 
and assume given a factorization
$(a \circ j)_! M_U \to M_\BA \to (a \circ j)_* M_U$ of 
the morphism
$a_* m : (a \circ j)_! M_U \to (a \circ j)_* M_U$ through an object $M_\BA$ of 
$\CHFDM$. 
Then $a_* \ujast M_U$ is canonically identified with a direct factor of $M_\BA$,
admitting a canonical complement.
\end{Thm}

\begin{Proof}
This is \cite[Cor.~2.5]{W4}.
\end{Proof}

The theory applies in particular when $a$ equals the structure morphism 
from $X$ to the base scheme $\BB$.

\begin{Def} \label{2Ia}
Assume that $X$ is proper over $\BA = \BB$. 
Denote by $a: X \to \BB$ the structure morphism of $X$.
Let $M_U \in \CHUM_{F,\partial w \ne 0,1} \, $. We call $a_* \ujast M_U$ the
\emph{intersection motive of $U$ relative to $X$ with coefficients in $M_U$}.
\end{Def}

Our terminology is motivated by one of the main results of \cite{W9}. It states
that on Chow motives of Abelian type, the cohomological \emph{Betti} 
\cite[D\'ef.~2.1]{Ay2} and \emph{$\ell$-adic realizations} 
\cite[Sect.~7.2, see in part.\ Rem.~7.2.25]{CD2} are compa\-tible with 
intermediate extensions (of motives,
and of perverse sheaves). For details, we refer to \cite[Thm.~7.2]{W9}. Since
the realizations are compatible with direct images, they therefore map
$a_* \ujast M_U$ to \emph{intersection cohomology} whenever $M_U$ is 
a Chow motive of Abelian type. \\   

It turns out that the comparison results from \cite{CD} allow to
connect the above to the notion of interior motive.

\begin{Prop} \label{2K}
Assume that $\BA = \BB = \Spec k$ for a field $k$
admitting strict resolution of singularities. 
Assume also that the structure
morphism $a: X \to \BB$ is proper, and that its restriction $a \circ j$ to
$U \subset X$ is smooth. Let $\pi: \C \to U$ 
proper and smooth (hence $\C$ is smooth over $k$). 
Assume $C$ to be quasi-projective over $k$.
Consider the Chow motive $\pi_* \one_{\C}$ over $U$.
Then the morphism
\[
a_* m : (a \circ j)_! \pi_* \one_{\C} \longto (a \circ j)_* \pi_* \one_{\C} 
\]
is canonically and $\ch_{d_\C} (\C \times_U \C)$-equivariantly 
isomorphic  ($d_\C:=$ the absolute dimension of $\C$) to the dual of the morphism
\[
u : \Mgm(\C) \longto \Mcgm(\C)  
\]
in $\DgM_F$ \cite[pp.~223--224]{V}. 
\end{Prop}

A few words of explanation are in order. First, by \cite[Cor.~16.1.6]{CD},
the triangulated category
$\DBcFkM$ is identified with the $F$-linear version of 
the triangulated category of \emph{geometrical motives} 
$\DgM_F$ \cite[p.~192]{V} 
(see \cite[Sect.~17.1.3]{A}). Second,
the duality in question is the functor mapping $N$ to
\[
N^* := \underline{Hom} \bigl( N, \one_{\Spec k} \bigl) \; .
\]
Third, equivariance under the \emph{Chow group} 
$\ch_{d_\C} (\C \times_U \C)$
refers to the following. 
According to
\cite[Cor.~14.2.14]{CD}, 
\[
\End_{CHM(U)_F} \bigl( \pi_* \one_{\C} \bigr)
= \ch_{d_\C} (\C \times_U \C) \otimes_\BZ F \; ,
\] 
meaning that the Chow group $\ch_{d_\C} (\C \times_U \C)$
acts on $\pi_* \one_{\C} \in \DBcFUM$. Hence the morphism $a_*m$ is
$\ch_{d_\C} (\C \times_U \C)$-equivariant. As for the action on 
$\Mgm(\C)$, on $\Mcgm(\C)$, and the equivariance of $u$, we refer to 
\cite[Thm.~5.23]{D} and \cite[Ex.~4.12, Ex.~7.15]{CD1}. 

\begin{Rem}
(a)~According to  \cite[Prop.~5.19, Cor.~6.14]{L2}, the identification 
\[
\End_{CHM(U)_F} \bigl( \pi_* \one_{\C} \bigr)
= \ch_{d_\C} (\C \times_U \C) \otimes_\BZ F \; ,
\] 
is compatible with composition. 
Thus, the action of the ring $\ch_{d_\C} (\C \times_U \C)$ on 
$\pi_* \one_{\C}$ is by correspondences in the classical sense. \\[0.1cm]
(b)~Note that since we assumed $C$ to be proper and smooth over $U$,
the Chow motive $\pi_* \one_{\C}$ is \emph{smooth} in the sense of \cite{L2}. 
According to \cite[Thm.~3.17]{F}, if $U$ is quasi-projective over $k$, then the comparision statement 
\[
\End_{CHM(U)_F} \bigl( \pi_* \one_{\C} \bigr)
= \ch_{d_\C} (\C \times_U \C) \otimes_\BZ F \; ,
\]
continues to hold if $C$ is assumed to be proper over $U$, and 
to remain quasi-projective and
smooth over $k$.
\end{Rem} 

\begin{Proofof}{Proposition~\ref{2K}}
The morphism $a_*m$ coincides with the value of the natural transformation
of functors
\[
(a \circ j \circ \pi)_! \longto (a \circ j \circ \pi)_*  
\]
on $\one_{\C}$, since both $\pi$ and $a$ are proper \cite[Thm.~2.2.14~(2)]{CD}.
Fix a projective and smooth compactification $\bar{\C}$ of $\C$ over $k$,
write $j'$ for the open immersion of $\C$ into $\bar{\C}$,
and $c$ for the structure morphism of $\bar{\C}$. 
Thus, $c \circ j' = a \circ j \circ \pi$ is the structure morphism of $\C$. 
Then $a_*m$ equals $c_* m'$, for the morphism
$m' : j'_! \one_{\C} \to j'_* \one_{\C}$ in $\DBcM(\bar{\C})_F$, 
and can be factorized as follows: 
\[
a_* m = c_* adj_1 \circ c_* adj_2 :
c_* j'_! \one_{\C} \stackrel{c_* adj_2}{\longto} 
c_* \one_{\bar{\C}} \stackrel{c_* adj_1}{\longto} 
c_* j'_* \one_{\C} \; ,
\]
where $adj_1 : \one_{\bar{\C}} \to j'_* (j')^* \one_{\bar{\C}} = j'_* \one_{\C}$ 
and 
$adj_2 : j'_! \one_{\C} = j'_! (j')^* \one_{\bar{\C}} \to \one_{\bar{\C}}$ 
are the adjunction maps.
Now $c_* adj_2$ and $c_* adj_1$ are related by duality: 
we have 
\[
c_* j'_! \one_{\C} = \bigl( c_* j'_* (c \circ j')^! \one_{\Spec k} \bigr)^* 
= \bigl( c_* j'_* (j')^* c^! \one_{\Spec k} \bigr)^* \; , 
\]
\[
c_* \one_{\bar{\C}} = \bigl( c_* c^! \one_{\Spec k} \bigr)^*
\; , \; 
\] 
and under these identifications, $c_* adj_2$ is dual to 
the morphism
\[
c_* adj_1^!: c_* c^! \one_{\Spec k} \longto 
c_* j'_* (j')^* c^! \one_{\Spec k} \; ,
\]
where $adj_1^!$ denotes the adjunction map 
$c^! \one_{\Spec k} \to j'_* (j')^* c^! \one_{\Spec k}$
\cite[Thm.~15.2.4]{CD}. Fix an isomorphism 
$\alpha: c^! \one_{\Spec k} \isoto \one_{\bar{\C}}(d_\C)[2d_\C]$; according to
\cite[Cor.~3.7]{W10}, such an isomorphism exists.
It is unique up to multiplication by global sections of
the constant sheaf $F^*$ on $C$. 
\emph{Via} $\alpha$, the morphism $c_* adj_1^!$ is identified
with 
\[
c_* adj_1(d_\C)[2d_\C] : 
c_* \one_{\bar{\C}}(d_\C)[2d_\C] \longto 
c_* j'_* (j')^* \one_{\bar{\C}}(d_\C)[2d_\C] \; ,
\]
and this identification does not
depend on the choice of $\alpha$.

To summarize the discussion so far: the morphism $a_* m$ equals the composition
of $c_* adj_1$, preceded by the dual of $c_* adj_1(d_\C)[2d_\C]$. 

As for the morphism $u$, observe that it, too, can be factorized:
\[
u = (j')^* \circ j_*' : \Mgm(\C) \stackrel{j_*'}{\longto} 
\Mgm (\bar{\C}) = \Mcgm (\bar{\C}) 
\stackrel{(j')^*}{\longto} \Mcgm(\C) \; ,
\]
where we denote by $j_*'$ and $(j')^*$ the morphisms induced by
the open immersion $j'$ on the level of $\Mgm$ and $\Mcgm$,
respectively \cite[pp.~223--224]{V}. 
According to \cite[Thm.~4.3.7~3]{V}, the dual of $(j')^*$
is identified with 
\[
j_*'(-d_\C)[-2d_\C] : 
\Mgm (\C)(-d_\C)[-2d_\C] \longto \Mgm (\bar{\C})(-d_\C)[-2d_\C] \; .
\]
To summarize: the morphism $u$ equals the composition of
$j_*'$, followed by the dual of $j_*'(-d_\C)[-2d_\C]$.

To relate $a_* m$ and $u$, observe that
\[
\Mgm(\C) = (c \circ j')_\sharp \one_{\C} = c_\sharp j_!' \one_{\C} \; , \;
\Mgm (\bar{\C}) = c_\sharp \one_{\bar{\C}} \; , \;
\]
and under these identifications, the morphism 
$j_*'$ equals $c_\sharp adj_3$, where 
$adj_3 : j_!' \one_{\C} = j_!' (j')^* \one_{\bar{\C}} \to \one_{\bar{\C}}$
is the adjunction \cite[Sect.~1.1.34, Sect.~11.1.2, Cor.~16.1.6]{CD}.
According to one of the projection formulae
\[
\underline{Hom}_V(f_\sharp M, N) = f_* \underline{Hom}_T(M,f^*N)
\]
for $f: T \to V$ smooth \cite[Sect.~1.1.33]{CD},
the morphism $j_*' = c_\sharp adj_3$ is dual to $c_* adj_1$.
Thus, the dual of $u$ equals $c_* adj_1$, preceded by
the dual of $c_* adj_1(d_\C)[2d_\C]$, \emph{i.e.}, 
it equals $a_* m$.

It remains to show that the identification of $a_*m$ and the dual of $u$ is
equivariant under $\ch_{d_\C} (\C \times_U \C)$. Given that $a_*m$
is the value at $\pi_* \one_{\C}$ of a natural transformation of functors,
and that $\ch_{d_\C} (\C \times_U \C) \otimes_\BZ F$ is identified with 
$\End_{CHM(U)_F} \bigl( \pi_* \one_{\C} \bigr)$, all one needs to establish is 
that under our identification, the action of $\ch_{d_\C} (\C \times_U \C)$
on $(a \circ j)_* \pi_* \one_{\C}$ coincides with the action of 
$\ch_{d_\C} (\C \times_U \C)$ on the dual of $\Mgm(\C)$, and likewise for
$(a \circ j)_! \pi_* \one_{\C} = \Mcgm(\C)^*$. As before, the second compatibility
is dual to the first, up to application of a twist by $d_\C$ and a shift by $2d_\C$.

As for the identification 
$(a \circ j)_* \pi_* \one_{\C} = c_* j'_* \one_{\C} = \Mgm(\C)^*$, it is compatible
with the action of \emph{finite correspondences} $\FZ \in c_{U}(\C,\C)$
by the very definition of the category of motivic complexes \cite[Def.~11.1.1]{CD}.
It remains to cite \cite[Lemma~5.18]{L2}: every class in
$\ch_{d_\C} (\C \times_U \C)$
can be represented by a cycle $\FZ$ belonging to 
$c_{U}(\C,\C)$.
\end{Proofof}
  
\begin{Cor} \label{2L}
Assume that $\BA = \BB = \Spec k$ for a field $k$
admitting strict resolution of singularities, that the structure
morphism $a: X \to \BB$ is proper, and that its restriction $a \circ j$ to
$U \subset X$ is smooth. Let $\pi: \C \to U$ 
proper and smooth, and assume that $C$ is quasi-projective over $k$.
Let $e \in \ch_{d_\C} (\C \times_U \C) \otimes_\BZ F$
an idempotent. 
Assume that the direct factor $(\pi_* \one_{\C})^e$ of the Chow motive 
$\pi_* \one_{\C} \in \CHUM_F$
lies in $\CHUM_{F,\partial w \ne 0,1} \, $. \\[0.1cm]
(a)~The $e$-part $\dMgm(\C)^e$ of the \emph{boundary motive} $\dMgm(\C)$ 
is without weights $-1$ and $0$.
In particular, $\C$ and $e$ satisfy assumption \cite[Asp.~4.2]{W4}, 
and therefore, the \emph{$e$-part of the interior motive of $\C$},
$\Gr_0 \Mgm(\C)^e$ is defined \cite[Def.~4.9]{W4}. \\[0.1cm]  
(b)~There is a canonical isomorphism 
\[
a_* \ujast (\pi_* \one_{\C})^e \isoto \bigl( \Gr_0 \Mgm(\C)^e \bigr)^* \; .
\] 
It is compatible with the factorizations
\[
(a \circ j)_! (\pi_* \one_{\C})^e \longto 
a_* \ujast (\pi_* \one_{\C})^e \longto
(a \circ j)_* (\pi_* \one_{\C})^e
\]
of $a_*m$ and 
\[
\Mgm(\C)^e \longto \Gr_0 \Mgm(\C)^e \longto \Mcgm(\C)^e
\]
of $u$ under the identification of Proposition~\ref{2K}.    
\end{Cor} 

\begin{Rem}
(a)~The hypothesis on strict resolution of singularities is (implicitly) used twice
(apart from the proof of Proposition~\ref{2K}). First, 
the results from \cite[Sect.~4]{W4} were formulated only for such fields.
This comes mainly from the fact
that at the time when \cite{W4} was written, 
the existence of the motivic weight structure on $\DgM_F$ was only established
under that hypothesis. Given the main results from \cite{Bo3}, 
one can dispose of that restriction on $k$ as far as the
weight structure is concerned (recall that our ring of coefficients $F$ is assumed to be 
a $\BQ$-algebra). 

Second, and more seriously, the hypothesis is used for the construction of 
the action of $\ch_{d_\C} (\C \times_U \C)$ on the boundary motive 
$\dMgm(\C)$ \cite[Thm.~2.2]{W6}, and hence for the very definition of $\dMgm(\C)^e$. 
It seems plausible that the hypothesis can be avoided using the main results from \cite[Sect.~5.3]{K},
in particular, localization for $\Mcgm \otimes_\BZ F$ \cite[Prop.~5.3.5]{K}, but we have not tried.

Given Corollary~\ref{2L}, the reader should obviously feel free to \emph{define}
the $e$-part of the interior motive of $\C$ as $(a_* \ujast (\pi_* \one_{\C})^e)^*$, 
in case the field $k$ does not
admit strict resolution of singularities. \\[0.1cm]
(b)~Recall from \cite[Def.~4.1~(a)]{W4} that there is
a ring $c_{1,2}(\C,\C)$
(of ``bi-finite correspondences'')
acting on the exact triangle 
\[
(\ast) \quad\quad
\dMgm(\C) \longto \Mgm(\C) \longto \Mcgm(\C) \longto \dMgm(\C)[1] \; .
\]
Denote by $\bar{c}_{1,2}(\C,\C)$ 
the quotient of $c_{1,2}(\C,\C)$ by the kernel of this action. 
The algebra $\bar{c}_{1,2}(\C,\C) \otimes_\BZ F$ is a canonical source
of idempotent endomorphisms of $(\ast)$, and it is for such choices
that assumption \cite[Asp.~4.2]{W4} was formulated. However, \cite[Asp.~4.2]{W4} 
admits an obvious
generalization to arbitrary idempotent endomorphisms $e$ of $(\ast)$. 
Similarly, \cite[Def.~4.9]{W4} and all
results from \cite[Sect.~4]{W4} remain valid in the present context,
up to modifications of the equivariance statement in \cite[Thm.~4.3]{W4} under
the centralizer of $e$ in $\bar{c}_{1,2}(\C,\C)$,
and of the explicit description of the effect of duality on $e$ in 
\cite[Prop.~4.15]{W4} (neither of which will be needed in the sequel). 
\end{Rem}

\begin{Proofof}{Corollary~\ref{2L}}
According to Proposition~\ref{2K}, 
\[
a_* m : (a \circ j)_! (\pi_* \one_{\C})^e \longto 
(a \circ j)_* (\pi_* \one_{\C})^e
\]
is identified with the dual of
\[
u : \Mgm(\C)^e \longto \Mcgm(\C)^e
\]
Any choice of cone of $a_*m$ is therefore isomorphic to 
the shift by $[1]$ of the dual of any choice of  
cone of $u$. But $a_* i_* i^* j_* (\pi_* \one_{\C})^e$ is a cone of $a_*m$,
while $\dMgm(\C)^e$ is the shift by $[-1]$ of a cone of $u$. Thus,
\[
\dMgm(\C)^e \cong \bigl( a_* i_* i^* j_* (\pi_* \one_{\C})^e \bigr)^* \; .
\]
Thanks to our additional assumption on $(\pi_* \one_{\C})^e$, and to
Proposition~\ref{2D}~(b), the motive $a_* i_* i^* j_* (\pi_* \one_{\C})^e$
is without weights $0$ and $1$. Part~(a) of our claim then follows from the
compatibility of the motivic weight structure with duality
\cite[Thm.~1.12]{W10}.

The same argument yields that $\Gr_0 (a \circ j)_* (\pi_* \one_{\C})^e$
and $\Gr_0 \Mgm(\C)^e$ are dual to each other. Part~(b) thus follows from
Theorem~\ref{2G}~(c).
\end{Proofof}

\begin{Rem}
An alternative proof could be given by showing that
the exact triangle
\[
a_* j_! (\pi_* \one_{\C}) \stackrel{a_* m}{\longto} a_* j_* (\pi_* \one_{\C}) 
\longto a_* i_* i^* j_* (\pi_* \one_{\C}) \longto a_* j_! (\pi_* \one_{\C}) [1]
\]
is $\ch_{d_\C} (\C \times_U \C)$-equivariantly 
isomorphic to the dual of the exact triangle
\[
\Mcgm(\C)[-1] \longto \dMgm(\C) \longto \Mgm(\C) \stackrel{u}{\longto}
\Mcgm(\C)\; .  
\]
To establish that latter result, one would apply techniques
similar to the ones used in the proofs of \cite[Thm.~2.2 and 2.5]{W6}.
\end{Rem}

\begin{Rem} \label{2M}
(a)~Assume that we are in the setting of Corollary~\ref{2L}. In particular, $\BA = \BB = \Spec k$,
and the (structure) morphism $a: X \to \BB$ is proper. From Corollary~\ref{2L}~(b),
and from \cite[Theorems~4.7 and 4.8]{W4}, it follows that the Chow motive
$a_* \ujast (\pi_* \one_{\C})^e$, \emph{i.e.}, the $e$-part of the intersection motive of $U$ 
relative to $X$ with coefficients in $\pi_* \one_{\C}$, 
realizes to give the $e$-part of \emph{interior cohomology}
of $C$. \\[0.1cm]
(b)~In fact, as is shown in \loccit, the statement from (a) follows from the more precise fact
that the values of the respective cohomological
(Hodge theoretic or $\ell$-adic) 
realization $(H^m \circ R)_{m \in \BZ}$ on the canonical morphisms 
\[
(a \circ j)_! (\pi_* \one_{\C})^e \longto 
a_* \ujast (\pi_* \one_{\C})^e 
\]
and 
\[
a_* \ujast (\pi_* \one_{\C})^e \longto
(a \circ j)_* (\pi_* \one_{\C})^e
\]
identify $H^m \circ R (a_* \ujast (\pi_* \one_{\C})^e)$ with the part of weight $m$ of
$H^m \circ R ((a \circ j)_! (\pi_* \one_{\C})^e)$, and of  
$H^m \circ R ((a \circ j)_* (\pi_* \one_{\C})^e)$, respectively. \\[0.1cm]
(c)~The statement from (b) can be shown without using Corollary~\ref{2L}, \emph{i.e.},
without any reference to the interior motive, 
and for arbitrary objects $M_U$ of $\CHUM_{F,\partial w \ne 0,1}$ instead of
$(\pi_* \one_{\C})^e$,
by formally imitating the proofs of \cite[Theorems~4.7 and 4.8]{W4}. 
The latter make essential use of the existence of weights on the level of realizations;
indeed, \cite[Theorems~4.7 and 4.8]{W4} should be seen as sheaf theoretic phenomena:
for any (Hodge theo\-retic or $\ell$-adic) sheaf
$N_U$ which is pure of weight $n$, and such that $i^*j_* N_U$ is without weights
$n$ and $n+1$, 
intersection and interior cohomology with coefficients in $N_U$ coincide. 
In particular, if $M_U \in \CHUM_{F,\partial w \ne 0,1}$ is of Abelian type, then
according to \cite[Thm.~7.2]{W9},
the natural map from  
intersection cohomology to cohomology 
identifies intersection cohomology and interior cohomology with coefficients
in the realization of $M_U$. \\[0.1cm]
(d)~For the $\ell$-adic realization, a relative version of statement~(c) holds,
provided that morphisms in the image of the realization are strict with respect
to the weight filtration: the
morphism $a: X \to \BA$ is still assumed to be proper, but $\BA$ may be different from
the base scheme, and the latter need not be a field. 
For a detailed study of the condition on strictness, we refer to \cite[Sect.~2]{Bo2};
note that it is satisfied in the situation we are about to study in Section~\ref{2}. \\[0.1cm]
(e)~An analogue of (d) should hold for the Hodge theoretic realization.
\end{Rem}

The following general result will be used in Section~\ref{3}.

\begin{Prop} \label{2O}
Let $g: W' \to W$ be a finite, \'etale morphism of schemes.
Then the direct image $g_*: \DBcFZoM \to \DBcFWM$ 
is \emph{weight conservative}, \emph{i.e.}, it detects weights.
More precisely, let $N' \in \DBcFZoM$, and $\alpha \le \beta$ two integers. 
\begin{enumerate}
\item[(a)] $N'$ lies in the heart $\CHFZoM$ if and only if 
$g_* N'$ lies in the heart $\CHFWM$.
\item[(b)] $N'$ lies in $\DBcM(W')_{F,w \le \alpha}$ if and only if 
$g_* N'$ lies in $\DBcM(W)_{F,w \le \alpha}$.  
\item[(c)] $N'$ lies in $\DBcM(W')_{F,w \ge \beta}$ if and only if 
$g_* N'$ lies in $\DBcM(W)_{F,w \ge \beta}$.
\item[(d)] $N'$ is without weights $\alpha,\alpha+1,\ldots,\beta$ if and only if 
$g_* N'$ is without weights $\alpha,\alpha+1,\ldots,\beta$.
\end{enumerate}
\end{Prop}

\begin{Proof}
Since $g$ is finite and \'etale, both $g^*$
and $g_*$ are weight exact \cite[Thm.~3.8~(ii'), (i), (i')]{H}.
In particular, the ``only if'' parts of statements (a)--(d) are true.

As for the ``if'' parts, note first that $g^* = g^!$ \cite[Thm.~2.4.50~(3), Def.~2.4.12~(2)]{CD}
and $g_! = g_*$ (since $g$ is proper). Therefore,
there are adjunction morphisms between $\id_{\DBcFZoM}$ and $g^* g_*$
in both directions.
Next, for any $N' \in \DBcFZoM$, the composition of the adjunctions
\[
N' \longto g^* g_* N' \longto N' 
\]
allows to
identify $N'$ with a direct factor of $g^* g_* N'$. 
Statements (a)--(c) thus follow from the fact
that the categories $\CHFZoM$, $\DBcM(W')_{F,w \le \alpha}$
and $\DBcM(W')_{F,w \ge \beta}$ are all pseudo-Abelian.
Statement~(d) is a consequence of functoriality of weight filtrations
avoiding weights $\alpha,\alpha+1,\ldots,\beta$ \cite[Prop.~1.7]{W4},
applied to an idempotent endomorphism of $g^* g_* N'$ cutting out $N'$.
\end{Proof}

\begin{Rem}
The analogue of Proposition~\ref{2O} holds for the inverse image 
$g^*$ under a finite, \'etale morphism $g$, with the same proof,
provided that $g$ is surjective.
This fact will not be needed in the sequel.
\end{Rem}

\begin{Cor} \label{2P}
Assume that $\BA = \BB = \Spec k$ for a field $k$
admitting strict resolution of singularities, that the structure
morphism $a: X \to \BB$ is proper, and that its restriction $a \circ j$ to
$U \subset X$ is smooth. Let $\pi: \C \to U$ 
proper and smooth, and assume that $C$ is quasi-projective over $k$.
Let $e \in \ch_{d_\C} (\C \times_U \C) \otimes_\BZ F$
an idempotent. 
Assume that the restriction $a \circ i$ of the structure morphism to $Z \subset X$
is finite. Let $\alpha \le \beta$ be two integers.
Then the following are equivalent.
\begin{enumerate}
\item[(1)] The motive $i^* j_* (\pi_* \one_{\C})^e \in \DBcFZM$
is without weights $\alpha,\alpha+1,\ldots,\beta$. 
\item[(2)] The motive $\dMgm(\C)^e \in \DBcFkM$
is without weights $- \beta,- (\beta - 1),\ldots,- \alpha$.
\end{enumerate}
In particular, the Chow motive $(\pi_* \one_{\C})^e$ lies in $\CHUM_{F,\partial w \ne 0,1}$
if and only if $\dMgm(\C)^e$ is without weights $-1$ and $0$.
\end{Cor}

\begin{Proof}
The field $k$ is perfect; therefore, the reduced scheme $Z_{red}$
underlying $Z$ is 
finite and \'etale over $\Spec k$. Given localization for the inclusion 
of $Z_{red}$ into $Z$ \cite[Sect.~2.3]{CD2}, we may thus assume that
$g := a \circ i$ is finite and \'etale. Given that
\[
\dMgm(\C)^e \cong \bigl( g_* i^* j_* (\pi_* \one_{\C})^e \bigr)^* \; .
\]
(see the proof of Corollary~\ref{2L}), our claim follows from 
Proposition~\ref{2O}.
\end{Proof}


\bigskip
%
%

\section{A criterion on absence of weights in the boundary}
\label{2}



We keep the geometrical situation of the preceding section: $X$ is a scheme, 
and $j: U \into X$ an open immersion with complement $i: Z \into X$. 
For a finite stratification by nilregular locally closed sub-schemes
$\Zp$ of $Z$, indexed by $\varphi \in \Phi$, recall the definition of
the category $\DBcFZPAbM$ of \emph{$\Phi$-constructible motives of Abelian type over $Z$} 
\cite[Def.~3.5~(b)]{W11}: it is the strict, full, dense, $F$-linear 
triangulated sub-category of $\DBcFZM$ 
generated by the images under $\pi_*$ of the objects of $DMT_\FS(S(\FS))_F$,
the \emph{category of $\FS$-constructible
Tate motives over $S(\FS)$} \cite[Def.~3.3]{W11}, where
\[
\pi: S(\FS) \longto Z = Z(\Phi)
\]
runs through the morphisms \emph{of Abelian type} with target
equal to $Z = Z(\Phi)$. 
According to \cite[Def.~3.5~(a)]{W11}, this means that $\FS$ is a finite
stratification of $S = S(\FS)$ by nilregular locally closed sub-schemes, that 
$\pi$ is a \emph{morphism of good stratifications} \cite[Def.~3.4]{W11},
that $\pi$ is proper, and that the following conditions are satisfied. 
\begin{enumerate}
\item[(1)] For any immersion 
$i_\tau: S_\tau \into \bSs$ of a stratum $S_\tau$ into
the closure $\bSs$ of a stratum $S_\sigma$, the functor $i_\tau^!$
maps $\one_{\bSs}$ to a Tate motive over $S_\tau$. 
\item[(2)] For all $\sigma \in \FS$
such that $\Ss$ is a stratum
of $\pi^{-1}(\Zp)$, the morphism
$\pis : \Ss \to \Zp$ can be factorized,
\[
\pis = \pios \circ \pits : \Ss \stackrel{\pits}{\longto} \Bs 
\stackrel{\pios}{\longto} \Zp \; ,
\]
such that the motive 
\[
\pi''_{\sigma,*} \one_{\Ss} \in \DBcFBsM
\]
belongs to the category $\DFTBsM$ of Tate motives over $\Bs \, $,
the morphism $\pios$ is proper and smooth, 
and its pull-back to any geometric
point of $\Zp$ lying over a generic point
is isomorphic to a finite disjoint union of Abelian varieties.
\end{enumerate}  

\begin{Def} 
An object $M \in \DBcFZM$ is a \emph{motive of Abelian type over $Z$} 
if it belongs to the sub-category $\DBcFZPAbM$, for a suitable
finite stratification $\Phi$ by nilregular locally closed sub-schemes of $Z$.
In this situation, we say that $\Phi$ is \emph{adapted to $M$}.
\end{Def}

Let us now fix a generic point $\Spec k$ of the base $\BB$. 
For any scheme $Y$,
denote by $R_{\ell,Y}$ the (generic) $\ell$-adic realization
\cite[Sect.~4]{W11}. Its target is the 
$F$-linear version $D^b_c(Y_k)_F$ 
of the bounded ``derived category'' $D^b_c(Y_k)$ \cite[Sect.~6]{E} of
constructible $\BQ_\ell$-sheaves on the fibre $Y_k$ of $Y$ over $\Spec k \into \BB$. 
According to \cite[Thm.~7.2.24]{CD2}, the $R_{\ell,\bullet}$ are compatible 
with the functors $f^* , f_* , f_! , f^!$. Furthermore, they are symmetric monoidal;
in particular, 
\[
R_{\ell,Y} (\one_Y) = 
\BQ_{\ell, Y_k} \; ,
\]
where $\BQ_{\ell, Y_k}$ denotes the $\ell$-adic structure sheaf on  $Y_k$. \\

\forget{
For a given $M_U \in \CHFUM$, consider the following assumption.

\begin{Ass} \label{Ass}
For all $n \in \BZ$, the perverse sheaf $H^n \! R_{\ell,Z}(i^*j_* M_U)$, 
given as the $n$-th perverse cohomology object of 
$R_{\ell,Z}(i^*j_* M_U)$, is without weights $n$ and $n+1$. 
\end{Ass}
}
The following is an immediate consequence of the main result from \cite{W11}. 

\begin{Thm} \label{2B}
Assume that $k$ is of characteristic zero.
Let $\ell$ be a prime. Let $N \in \DBcFZPAbM \, $.
Assume that the generic points of all strata $\Zp$, $\varphi \in \Phi$, 
lie over $\Spec k \into \BB$. For $\varphi \in \Phi$, denote by $\ip$ the
immersion of $\Zp$ into $Z$. \\[0.1cm]
(a)~Let $\alpha \in \BZ$. Then $N$ 
lies in $DM_{\text{\cyrb},c,\Phi}^{Ab}(Z)_{F,w \le \alpha}$ if and only if 
for all $n \in \BZ$, and all $\varphi \in \Phi$, 
the perverse cohomology sheaf 
\[
H^n i_\varphi^* R_{\ell,Z}(N)
\]
is of weights $\le n + \alpha$. \\[0.1cm]
(b)~Let $\beta \in \BZ$. Then $N$
lies in $DM_{\text{\cyrb},c,\Phi}^{Ab}(Z)_{F,w \ge \beta}$ if and only if 
for all $n \in \BZ$, and all $\varphi \in \Phi$, 
the perverse cohomology sheaf
\[
H^n i_\varphi^! R_{\ell,Z}(N)
\]
is of weights $\ge n + \beta$. 
\end{Thm}

\begin{Rem}
(a)~The conditions on the generic points of the strata 
$\Zp$ are empty when $\BB$ is itself the spectrum
of a field of characteristic zero. \\[0.1cm]
(b)~Recall from \cite[Prop.~2.1.2~1]{Bo} that any additive functor
$\CH$ from a triangulated category $\CC$ carrying a weight structure $w$,
to an Abelian category $\FA$ admits a canonical \emph{weight filtration}
by sub-functors
\[
\ldots \subset W_n \CH \subset W_{n+1} \CH \subset \ldots \subset \CH \; .
\]
For any $m \in \BZ$, one defines 
\[
\CH^m : \CC \longto \FA \; , \; M \longto \CH \bigl( X[m] \bigr) \; ;
\]
according to the usual convention, the weight filtration of $\CH^m(M)$ \emph{equals} the
weight filtration of $\CH(X[m])$, \emph{i.e.}, it differs by \emph{d\'ecalage}
from the intrinsic weight filtration of the covariant additive functor $\CH^m$.

If $k$ is finitely generated over $\BQ$, then there is an instrinsic notion
of weights on those perverse sheaves on $Y_k$, which are in the
image of the cohomological realization \cite[Prop.~2.5.1~(II)]{Bo2}.

In general, the weights of $H^* \! R_{\ell,Y}(M)$ are by definition
those induced by the weight filtration of the functor
$H^* \! R_{\ell,Y}$ 
(these coincide with the
above when $k$ is finitely generated over $\BQ$). 
\end{Rem}

\medskip

\begin{Proofof}{Theorem~\ref{2B}}
The motive $N$ belongs to 
$DM_{\text{\cyrb},c,\Phi}^{Ab}(Z)_{F,w \le \alpha}$ if and only if 
for all $\varphi \in \Phi$, 
\[
i_\varphi^* N \in DM_{\text{\cyrb},c,\Phi}^{Ab}(\Zp)_{F,w \le \alpha} \; .
\]
According to \cite[Thm.~4.4~(b)]{W11}, the latter condition is equivalent to
\[
H^n R_{\ell,\Zp} \bigl( i_\varphi^* N \bigr) \quad \text{is of weights} \quad \le n + \alpha \; ,
\]
for all $n \in \BZ$. But thanks to 
the compatibility of $R_{\ell,\bullet}$ with $i_\varphi^*$,
we have 
\[
R_{\ell,\Zp} \bigl( i_\varphi^* N \bigr) = i_\varphi^* R_{\ell,Z}(N) \; .
\]
This proves part~(a) of our claim. Dualizing, we obtain the proof of part~(b). 
\end{Proofof}

Together with one of the main compatibility results from \cite{W9}, we obtain
the following.

\begin{Thm} \label{2C}
Assume that $k$ is of characteristic zero.
Let $\ell$ be a prime. 
Let $M_U \in \CHFUM$, such that $R_{\ell,U}(M_U)$ is concentrated in a single perverse
degree, and such that $i^*j_* M_U \in \DBcFZM$ is of Abelian type. Let
$\Phi$ a stratification of $Z$ adapted to $i^*j_* M_U$.
Assume that the generic points of all $\Zp$, $\varphi \in \Phi$, lie over $\Spec k \into \BB$. 
For $\varphi \in \Phi$, denote by $\ip$ the
immersion of $\Zp$ into $Z$. \\[0.1cm]
(a)~The motive $M_U$ belongs to $\CHUM_{F,\partial w \ne 0,1}$ 
if and only if for all $n \in \BZ$, and all $\varphi \in \Phi$, the following 
conditions hold:
the perverse cohomology sheaf
\[
H^n i_\varphi^* i^* \ujast R_{\ell,U}(M_U)
\]
is of weights $\le n- 1$, and  
\[
H^n i_\varphi^! i^! \ujast R_{\ell,U}(M_U)
\]
is of weights $\ge n + 1$. 
In particular, the intermediate extension $\ujast M_U$ is then defined
up to unique isomorphism, as a Chow motive over $X$. \\[0.1cm]
(b)~Let $\alpha \le 0$ and $\beta \ge 1$ two integers.
The motive $i^*j_*M_U$ is without weights $\alpha , \alpha+1 , \ldots , \beta$
if and only if 
for all $n \in \BZ$, and all $\varphi \in \Phi$, the following conditions hold:
the perverse cohomology sheaf
\[
H^n i_\varphi^* i^* \ujast R_{\ell,U}(M_U)
\]
is of weights $\le n + \alpha - 1$, and  
\[
H^n i_\varphi^! i^! \ujast R_{\ell,U}(M_U)
\]
is of weights $\ge n + \beta$. \\[0.1cm]
(c)~Let $\alpha \le 0$ and $\beta \ge 1$ two integers.
Assume that for all $n \in \BZ$, and all $\varphi \in \Phi$, the following 
conditions hold:
the perverse cohomology sheaf
\[
H^n i_\varphi^* i^* \ujast R_{\ell,U}(M_U)
\]
is of weights $\le n + \alpha - 1$, and  
\[
H^n i_\varphi^! i^! \ujast R_{\ell,U}(M_U)
\]
is of weights $\ge n + \beta$. Then $a_* j_! M_U$ is without weights
$\alpha-1, \alpha, \ldots, -1$, 
and $a_* j_* M_U$ is without weights
$1, 2, \ldots, \beta$, for any proper morphism $a: X \to \BA$.
\end{Thm}

By slight abuse of notation, we write $\ujast R_{\ell,U}(M_U)$ for
\[
\bigl( \ujast (R_{\ell,U}(M_U)[s]) \bigr) [-s] \; ,
\]
if $R_{\ell,U}(M_U)$ is concentrated in perverse degree $s$.

\medskip

\begin{Proofof}{Theorem~\ref{2C}}
Part~(a) follows from part~(b)  (take $\alpha = 0$ and $\beta = 1$),
and from Definitions~\ref{0A}~(a) and \ref{0C}, while part~(c) is implied by (b) and 
Theorem~\ref{2G}~(a) and (b).

As for part~(b), 
note that according to \cite[Thm.~3.10]{W11}, the heart of the weight structure
of $\DBcFZPAbM$ is semi-primary. Therefore, Theorem~\ref{0E} can be applied; 
the motive $i^*j_*M_U$ is thus without weights $\alpha , \alpha+1 , \ldots , \beta$
if and only if 
\[
i^* \ujast M_U \in \CC(Z)_{w \le \alpha - 1} \quad \text{and} \quad
i^! \ujast M_U \in \CC(Z)_{w \ge \beta} \; .
\]
By Theorem~\ref{2B}, this is in turn equivalent to the following: 
for all $n \in \BZ$, and all $\varphi \in \Phi$,
\[
H^n i_\varphi^* R_{\ell,Z} \bigl( i^* \ujast M_U \bigr)
\]
is of weights $\le n + \alpha - 1$, and 
\[
H^n i_\varphi^! R_{\ell,Z} \bigl( i^* \ujast M_U \bigr)
\]
is of weights $\ge n + \beta$. Thanks to compatibility of 
$R_{\ell,\bullet}$ with $i^*$ and $i^!$,
we have 
\[
i^* \circ R_{\ell,X} = R_{\ell,Z} \circ i^* \quad \text{and} \quad
R_{\ell,Z} \circ i^! = i^! \circ R_{\ell,X} \; .
\]
The compatibility of $R_{\ell,\bullet}$ with $\ujast$ is the content of \cite[Thm.~7.2~(b)]{W9}.
\end{Proofof}

\begin{Rem} \label{2CR}
(a)~Given the full, triangulated sub-category $\DBcFZPAbM$ of $\DBcFZM \,$, there is
first a maximal choice of full, triangulated sub-category 
$\CD(U)$ of $\DBcFUM \,$, which glues with $\DBcFZPAbM \,$,
to give a full, triangulated sub-category of $\DBcFXM \,$, namely the full sub-category of objects
$M_U$ satisfying $i^*j_* M_U \in \DBcFZPAbM$ (cmp.\ \cite[Prop~4.1]{W9}).
Inside $\CD(U)$, we then find the maximal choice of full, triangulated sub-category 
$\CC(U)$ of $\DBcFUM \,$, which glues with $\DBcFZPAbM \,$, 
to give a full, triangulated sub-category $\CC(X)$ of $\DBcFXM \,$,
and which in addition inherits a weight
structure from the motivic weight structure on $\DBcFUM \,$, namely the full triangulated
sub-category generated by objects
$M_U$ of $\CHFUM$ satisfying $i^*j_* M_U \in \DBcFZPAbM \,$. 
The theory from \cite[Sect.~2]{W9} and from the present
Section~\ref{0} can thus be applied to 
the triplet of categories $\CC(U)$, 
$\CC(X) \subset \DBcFXM \,$, and $\CC(Z) = \DBcFZPAbM \,$. \\[0.1cm]
(b)~In particular, if an object $M_U$ of $\CHFUM$ is such that $i^*j_* M_U \in \DBcFZPAbM \,$,
\emph{i.e.}, $M_U$ belongs to $\CC(U)_{w=0}$, where $\CC(U)$ is as in (a), then 
($\ujast M_U$ exists, and) $i^* \ujast M_U$ and $i^! \ujast M_U$ belong to 
$\DBcFZPAbM$. This fact was implicitly used in the proof of Theorem~\ref{2C}. \\[0.1cm]
(c)~\emph{A priori}, the application of \cite[Thm.~7.2]{W9} necessitates the validity
of \cite[Asp.~7.1]{W9}: the motive $M_U$ belongs to the triangulated sub-category
$\pi_* DMT_{\FS_U}(S(\FS_U))_{F,w = 0}^\natural$ of $\CHFUM$, for an extension
of $\pi: S(\FS) \to Z$ to $X$. While \cite[Asp.~7.1~(b), (c)]{W9}
belong to the hypotheses of Theorem~\ref{2C}, \cite[Asp.~7.1~(a)]{W9}
is replaced by the condition that $M_U$ belong to the ca\-te\-gory $\CC(U)$ from~(a).
The proof of \cite[Thm.~7.2~(b)]{W9} carries over to this more general context
without any modification.
\end{Rem}

\begin{Cor} \label{2Ca}
Assume that $k$ is of characteristic zero.
Let $\ell$ be a prime. 
Let $M_U \in \CHFUM$, such that $R_{\ell,U}(M_U)$ is concentrated in a single perverse
degree $s$, and such that $R_{\ell,U}(M_U)$
is auto-dual up to a shift and a twist: 
\[
\BD_{\ell,U} \bigl( R_{\ell,U}(M_U) \bigr) \cong R_{\ell,U}(M_U)(s)[2s] 
\] 
($\BD_{\ell,U} =$ $\ell$-adic local duality on $U$).
Assume in addition that 
$i^*j_* M_U$ is of Abelian type. Let
$\Phi$ a stratification of $Z$ adapted to $i^*j_* M_U$.
Assume that the generic points of all $\Zp$, $\varphi \in \Phi$, lie over $\Spec k \into \BB$. 
For $\varphi \in \Phi$, denote by $\ip$ the
immersion of $\Zp$ into $Z$. \\[0.1cm]
(a)~The motive $M_U$ belongs to $\CHUM_{F,\partial w \ne 0,1}$ 
if and only if for all $n \in \BZ$, and all $\varphi \in \Phi$, the following holds:
the perverse cohomology sheaf
\[
H^n i_\varphi^* i^* \ujast R_{\ell,U}(M_U)
\]
is of weights $\le n- 1$. \\[0.1cm]
(b)~Let $\beta \ge 1$ be an integer. The following are equivalent.
\begin{enumerate}
\item[(b1)] The motive $i^*j_*M_U$ is without weights $-\beta+1 , -\beta+2 , \ldots , \beta$.
\item[(b2)] For all $n \in \BZ$, and all $\varphi \in \Phi$, 
\[
H^n i_\varphi^* i^* \ujast R_{\ell,U}(M_U)
\]
is of weights $\le n - \beta$. 
\item[(b3)] For all $n \in \BZ$, and all $\varphi \in \Phi$, 
\[
H^n i_\varphi^! i^! \ujast R_{\ell,U}(M_U)
\]
is of weights $\ge n + \beta$. 
\end{enumerate}
(c)~Let $\beta \ge 1$ be an integer, and assume that 
one of the equivalent conditions (b1), (b2), (b3) is satisfied.
Then $a_* j_! M_U$ is without weights
$-\beta, -\beta +1 , \ldots, -1$, 
and $a_* j_* M_U$ is without weights
$1, 2, \ldots, \beta$, for any proper morphism $a: X \to \BA$.
\end{Cor}

\begin{Proof}
Part~(a) is a special case of part~(b) (take $\beta = 1$, and use the equivalence
(b1)~$\Leftrightarrow$~(b2)). Similarly, part~(c) follows from 
(b2)~$\Leftrightarrow$~(b3), and from Theorem~\ref{2C}~(c) (with $\alpha = -\beta+1$).

As for part~(b), observe that
for all $n \in \BZ$, and all $\varphi \in \Phi$, 
\[
H^n i_\varphi^! i^! \ujast R_{\ell,U}(M_U)
\]
is dual to 
\[
H^{-n} i_\varphi^* i^* \ujast \BD_{\ell,U} \bigl( R_{\ell,U}(M_U) \bigr) \; .
\]
By assumption, weight $w$ occurs in 
\[
H^{-n} i_\varphi^* i^* \ujast \BD_{\ell,U} \bigl( R_{\ell,U}(M_U) \bigr) \cong
H^{-n} i_\varphi^* i^* \ujast \bigl( R_{\ell,U}(M_U)(s)[2s]  \bigr)
\]
if and only if weight $2s+w$ occurs in 
\[
H^{2s-n} i_\varphi^* i^* \ujast R_{\ell,U}(M_U) \; .
\]
Therefore, conditions~(b2) and (b3) are equivalent to each other. 
According to Theorem~\ref{2C}~(b)
(with $\alpha = -\beta+1$), each of them is thus equivalent to (b1). 
\end{Proof}

\begin{Rem} \label{2Cb}
Applying the
variant of the theory from Section~\ref{0} sketched in Remark~\ref{0F},
we see that in Corollary~\ref{2Ca}~(b), conditions~(b1)--(b3) are also equivalent to
\begin{enumerate}
\item[(b4)] The motive $i^*j_*M_U$ is without weights $-\beta+1 , -\beta+2 , \ldots , 0$.
\item[(b5)] The motive $i^*j_*M_U$ is without weights $1 , 2 , \ldots , \beta$.
\end{enumerate}
Clearly condition~(b1) implies both (b4) and (b5). We claim that
(b4) implies (b2), and that (b5) implies (b3). Indeed, according to
Remark~\ref{0F}~(b), condition~(b4) implies
\[
i^* \ujast M_U \in \CC(Z)_{w \le -\beta} \; ,
\]
while condition~(b5) implies
\[
i^! \ujast M_U \in \CC(Z)_{w \ge \beta} \; .
\]
Now apply Theorem~\ref{2B}~(a) and (b).
\end{Rem}

Together with the comparison result from Section~\ref{1}, we get the following.

\begin{Thm} \label{2N}
Assume that $\BB = \Spec k$ for a field $k$ of characteristic zero,
that the structure
morphism $a: X \to \BB$ is proper, and that its restriction $a \circ j$ to
$U \subset X$ is smooth. Let $\pi: \C \to U$ 
proper and smooth, and assume that $C$ is quasi-projective over $k$. 
Let $e \in \ch_{d_\C} (\C \times_U \C) \otimes_\BZ F$
an idempotent. Let $\ell$ be a prime. 
Assume that $R_{\ell,U}(\pi_* \one_{\C})^e$ is concentrated in a single perverse
degree $s$, and that $R_{\ell,U}(\pi_* \one_{\C})^e$
is auto-dual up to a shift and a twist: 
\[
\BD_{\ell,U} \bigl( R_{\ell,U}(\pi_* \one_{\C})^e \bigr) 
\cong R_{\ell,U}(\pi_* \one_{\C})^e(s)[2s] \; .
\] 
Assume in addition that the motive
$i^*j_* (\pi_* \one_{\C})^e$ is of Abelian type. Let
$\Phi$ a stratification of $Z$ adapted to $i^*j_* (\pi_* \one_{\C})^e$.
For $\varphi \in \Phi$, denote by $\ip$ the
immersion of $\Zp$ into $Z$. \\[0.1cm]
(a)~If for all $n \in \BZ$, and all $\varphi \in \Phi$, 
the perverse cohomology sheaf
\[
H^n i_\varphi^* i^* \ujast R_{\ell,U}(\pi_* \one_{\C})^e
\]
is of weights $\le n- 1$, then $\dMgm(\C)^e$ 
is without weights $-1$ and $0$. \\[0.1cm]  
(b)~If for all $n \in \BZ$, and all $\varphi \in \Phi$, 
the perverse cohomology sheaf
\[
H^n i_\varphi^* i^* \ujast R_{\ell,U}(\pi_* \one_{\C})^e
\]
is of weights $\le n- 1$, then there is a canonical isomorphism 
\[
a_* \ujast (\pi_* \one_{\C})^e \isoto \bigl( \Gr_0 \Mgm(\C)^e \bigr)^* \; .
\] 
\end{Thm}

\begin{Proof}
Combine Corollaries~\ref{2Ca} and \ref{2L}.
\end{Proof}

\begin{Rem}
The condition on auto-duality of $R_{\ell,U}(\pi_* \one_{\C})^e$ is satisfied
if the cycle $e$ and its transposition ${}^t e$ have identical images under
$R_{\ell,C}$.
\end{Rem}

\forget{
Throughout this section, we fix the following geometrical situation. The map
$\pi: S(\FS) \to Y(\Phi)$ is a 
\emph{morphism of good stratifications} in the sense of \cite[Sect.~4]{W9}:
the stratifications of $S(\FS)$ and $Y(\Phi)$ are such that the closure of any
stratum is a union of strata, and the same is true for the pre-image of any
stratum under $\pi$.
Furthermore, we assume that $\pi$ is proper (as a morphism of schemes).
We also fix an open subset $\Phi_U$ of $\Phi$, meaning that it gives rise to
a good stratification $Y(\Phi_U)$, and that $j : Y(\Phi_U) \into Y(\Phi)$
is an open immersion. Denote
by $\Phi_Z$ the complement of $\Phi_U$ in $\Phi$, and by $i$
the closed immersion complementary to $j$. Write $\FS_U := \pi^{-1} \Phi_U$ and
$\FS_Z := \pi^{-1} \Phi_Z$. Finally, we assume that
for all $\sigma \in \FS$, the closures $\bSs$ of strata $\Ss$ are 
nilregular, and that the induced  \\

According to \cite[Prop.~3.9]{W11}, the morphism $\pi_Z$ satisfies
the assumptions of \cite[Main Thm.~5.4]{W9}; therefore, $\pi$
satisfies \cite[Asp.~5.6]{W9}.
According to \cite[Cor.~5.7]{W9}, the hypotheses of
\cite[Thm.~2.9~(a)]{W9} are fulfilled, 
and the properties listed in \cite[Summ.~2.12]{W9} are satisfied by
the full, dense, triangulated sub-categories 
$\CC(X) = \pi_* DMT_\FS(S(\FS))_F^\natural \subset \DBcM(S(\FS))_F$, 
$\CC(U) = \pi_* DMT_\FS(S(\FS_U))_F^\natural \subset \DBcM(S(\FS_U))_F$, and
$\CC(Z) = \pi_{Z,*} DMT_\FS(S(\FS_Z))_F^\natural \subset \DBcM(S(\FS_Z))_F$.
Thus, the general theory of intermediate extensions 
from \cite[Sect.~2]{W9} is at our disposal.

\begin{Thm} \label{2A}
(a)~The restriction of the intermediate extension induces a 
fully faithful functor
\[
\pi_* DMT_\FS(S(\FS_U))_{F,w=0,\partial w \ne 0,1}^\natural \longinto
\pi_* DMT_\FS(S(\FS))_{F,w=0}^\natural \; ,
\]
denoted by the same symbol $\ujast$. It is the unique functor fitting into
a commutative diagram
\[
\vcenter{\xymatrix@R-10pt{
        \pi_* DMT_\FS(S(\FS_U))_{F,w=0,\partial w \ne 0,1}^\natural 
               \ar@{_{ (}->}[d] \ar[r]^-{\ujast} &
        \pi_* DMT_\FS(S(\FS))_{F,w=0}^\natural
               \ar@{->>}[d] \\
        \pi_* DMT_\FS(S(\FS_U))_{F,w=0}^\natural \ar@{^{ (}->}[r]^-{\ujast} &
        \pi_* DMT_\FS(S(\FS))_{F,w=0}^{\natural,u}         
\\}}
\] 
(b)~The essential image of 
$\pi_* DMT_\FS(S(\FS_U))_{F,w=0,\partial w \ne 0,1}^\natural$ under $\ujast$
is the full
sub-category $\pi_* DMT_\FS(S(\FS))_{F,w=0,i^*w \le -1,i^!w \ge 1}^\natural$ 
of $\pi_* DMT_\FS(S(\FS))_{F,w=0}^\natural$ consisting of objects $M$ 
such that both $i^* M$ and $i^! M$ are without weight $0$. The induced functor
\[
\ujast : \pi_* DMT_\FS(S(\FS_U))_{F,w=0,\partial w \ne 0,1}^\natural
\longto \pi_* DMT_\FS(S(\FS))_{F,w=0,i^*w \le -1,i^!w \ge 1}^\natural
\]
is an equivalence of categories, with inverse equal to $j^*$. \\[0.1cm]
(c)~Any object $M$ of $\pi_* DMT_\FS(S(\FS))_{F,w=0}^\natural$
is isomorphic to $\ujast M_U \oplus i_* N$,
for suitable $M_U \in \pi_* DMT_\FS(S(\FS_U))_{F,w=0}^\natural$
and $N \in \pi_{Z,*} DMT_\FS(S(\FS_Z))_{F,w=0}^\natural \, $.
The motive $M$ belongs to 
$\pi_* DMT_\FS(S(\FS))_{F,w=0,i^*w \le -1,i^!w \ge 1}^\natural$
if and only if on the one hand, $M_U$ belongs to 
$\pi_* DMT_\FS(S(\FS_U))_{F,w=0,\partial w \ne 0,1}^\natural$
and on the other, $N = 0$. 
\end{Thm}

\begin{Proof}
For part~(a) of the claim, see Remark~\ref{0D}~(d). 

For $M \in \pi_* DMT_\FS(S(\FS))_{F,w=0}^\natural \, $, 
the motives $i^* M$ and $i^! M$ are \emph{a priori}
of weights $\le 0$ and $\ge 0$, respectively.
To say that they are without weight $0$ is therefore equivalent to the
conditions
\[
i^* M \in \pi_{Z,*} DMT_\FS(S(\FS_Z))_{F,w \le -1}^\natural 
\quad \text{and} \quad
i^! M \in \pi_{Z,*} DMT_\FS(S(\FS_Z))_{F,w \ge 1}^\natural \; .
\]  
Part~(b) of our claim thus follows from the definitions,
and from Theorem~\ref{0B}.

Part~(c) is \cite[Summ.~2.12~(b)]{W9} and Remark~\ref{0D}~(a).
\end{Proof}

Together with the results from Section~\ref{1}, we get the following.
 
\begin{Cor} \label{2C}
Assume that $\BB = \Spec k$ for a field $k$
of characteristic zero, that the structure
morphism $d: X \to \BB$ is proper, and that its restriction to $U$ is smooth.
Let $\pi: \C \to U$ proper and smooth, and 
$e \in \ch_{d_\C} (\C \times_U \C) \otimes_\BZ F$ an idempotent. 
Let $\alpha \le 0$ and $\beta \ge 1$ two integers.
Assume that $i^*j_* (\pi_* \one_\C)^e \in \DBcFZM$ is of Abelian type,
and that for all $n \in \BZ$, the $e$-part 
\[
H^n \! R_{\ell,Z} \bigl( i^*j_* (\pi_* \one_\C) \bigr)^e
\]
of the perverse cohomology sheaf
$H^n \! R_{\ell,Z} ( i^*j_* (\pi_* \one_\C))$
is without weights $n + \alpha,n + \alpha+1,\ldots,n + \beta$. \\[0.1cm]
(a)~Assumption~\ref{Ass} is satisfied, for $M_U = (\pi_* \one_\C)^e$. 
Therefore, \cite[Asp.~4.2]{W4} is satisfied, for $C$ and $e$, 
the $e$-part of the interior motive of $\C$,
$\Gr_0 \Mgm(\C)^e$ exists, and 
\[
d_* \ujast (\pi_* \one_\C)^e \isoto \bigl( \Gr_0 \Mgm(\C)^e \bigr)^* \; .
\]
(b)~The motive $\Mcgm(\C)^e$ is without weights $1, \ldots, -\alpha, -\alpha+1$ and 
$\Mgm(\C)^e$ is without weights $-\beta, \ldots, -2, -1$. 
\end{Cor}

\begin{Proof}
(a): Apply Theorem~\ref{2B}~(b) and Corollary~\ref{2L}.

(b): Apply Theorem~\ref{2B}~(c), Proposition~\ref{2K},
and compatibility of the motivic weight structure with duality
\cite[Thm.~1.12]{W10}.
\end{Proof} 
}


\bigskip

%
%

\section{Examples: the boundary of certain Shimu\-ra varie\-ties}
\label{3}



The common features of the examples to be reviewed in the present section
are the following. 
The open immersion $j$ equals the inclusion
of a (pure) \emph{Shimura variety} of $PEL$-type 
$M^K$, whose \emph{level} $K$ is \emph{neat}, 
into its \emph{Baily--Borel 
compactification} $(M^K)^*$. The complement $i$ thus equals the closed
immersion of the boun\-dary $\partial (M^K)^*$ of $(M^K)^*$. 
The variety $M^K$ is associated to (pure) \emph{Shimura data} $(G,\CH)$;
in particular, $G$ is a connected reductive group over $\BQ \,$.
The finite stratification $\Phi$ of $\partial (M^K)^*$ is indexed by
the $G(\BQ)$-conjugation classes of \emph{rational boundary components}
of $(G,\CH)$. Each stratum is a finite disjoint union of 
locally closed sub-varieties of $\partial (M^K)^*$, each of which is
a quotient by the action of a finite group of a (pure) Shimura variety
associated to Shimura data, which are ``smaller'' than $(G,\CH)$. 
All strata are nilregular. The category
$\DBcFpPAbM$ of $\Phi$-constructible motives of Abelian type over $\partial (M^K)^*$
is therefore defined. \\

As for the source of the relative Chow motives in $\CHFMM$,
note that the Shimura data $(G,\CH)$ being of $PEL$-type, there is on the one hand
a canonical faithful representation $V$ of $G$ (the latter being defined as the group
of endomorphisms of $V$ commuting with a certain semi-simple algebra, and
respecting, up to scalars, a certain anti-symmetric bilinear form on $V$). 
On the other hand, given the modular interpretation of $M^K$, there
is a universal Abelian scheme $\pi: \B \to M^K$. 
Denote by $\pi_*^m \one_{\B}$, $m \ge 0$, the 
$m$-th \emph{Chow-K\"unneth component} of the Chow motive 
$\pi_* \one_{\B}$ over $M^K$  \cite[Thm.~3.1]{DM}. 

\begin{Thm}[{\cite[Thm.~8.6]{Anc}}] \label{3A}
There is an $F$-linear tensor functor
\[
\widetilde{\mu} : \Rep (G_F) \longto CHM^s (M^K)_F
\]
from the category $\Rep (G_F)$ of algebraic
representations of $G$ in $F$-modules of finite type
to the full sub-category $CHM^s (M^K)_F$ of $\CHFMM$
of smooth Chow motives over $M^K$.  
It has the following properties.
\begin{enumerate}
\item[(a)] The composition of $\widetilde{\mu}$
with the cohomological
Hodge theoretic reali\-zation is isomorphic to the \emph{canonical construction} 
functor $\mu_{\bf H}$ 
(e.g.\ \cite[Thm.~2.2]{W1}) to the category
of admissible graded-polarizable variations of Hodge structure on $M^K_\BC$.
\item[(b)] The composition of $\widetilde{\mu}$
with the cohomological
$\ell$-adic reali\-zation is isomorphic to the \emph{canonical construction} functor $\mu_\ell$ 
(e.g.\ \cite[Chap.~4]{W1}) to the category
of lisse $\ell$-adic sheaves on $M^K$.
\item[(c)] The functor $\widetilde{\mu}$ commutes with Tate twists
in the following sense:
for any $W \in \Rep (G_F)$ and $n \in \BZ$, we have
\[
\widetilde{\mu} \bigl( W(n) \bigr) = \widetilde{\mu} (W)(n)[2n] \; .
\]
\item[(d)] The functor $\widetilde{\mu}$
maps the representation $V$
to the Chow motive $\pi_*^1 \one_{\B}(1)[2]$ over $M^K$.
\end{enumerate} 
\end{Thm}

\begin{Proof}
Parts~(a), (c) and (d) are identical to \cite[Thm.~8.6]{Anc};
as for (d), note that the anti-symmetric bilinear form implicit in the $PEL$-data
induces an isomorphism between the dual of $V$ and $V(-1)$.

As for part~(b), repeat the proof of \loccit , observing that the $\ell$-adic analogue
of \cite[Prop.~8.5]{Anc} holds (the base change to $\BQ_\ell$ of the
sub-group $G_1$ of $G$ coincides with the Lefschetz group).   
\end{Proof}

Additional common features of Shimura varieties are
that the base scheme $\BB$ equals $\Spec E$, for a number field $E$
called the \emph{reflex field} of $(G, \CH)$,
and that the ring of coefficients
$F$ is equipped with a canonical
structure of $F'$-algebra, for a number field $F'$ over which
$G$ is split.

\begin{Def} \label{3B}
Fix a maximal split torus $T$ of $G_{F'}$, and a dominant character $\ua$ of $T$. \\[0.1cm]
(a)~Denote by $V_{\ua} \in \Rep (G_F)$ the irreducible representation
of highest weight $\ua$. \\[0.1cm]
(b)~Define ${}^{\ua} \CV \in \CHsMM_F \subset \CHFMM$ as 
\[
{}^{\ua} \CV := \widetilde{\mu}(V_{\ua}) \; .
\]
\end{Def}

\begin{Ex} \label{3C}
Our first example concerns \emph{modular curves}.
The reductive group $G$ equals $GL_{2,\BQ} \,$,
the reflex field equals $\BQ \,$, meaning that 
the base scheme $\BB$ equals $\Spec \BQ \,$, and $F = F' = \BQ \,$.
The dominant character $\ua$ is identified with a pair of integers $(k,r)$,
with $k \ge 0$ and $r \equiv k \mod 2$: choosing $T \subset GL_{2,\BQ}$ to be equal to 
the sub-group of diagonal matrices, we have
\[
\ua : T \longto \Gm \; , \; \diag(a,a^{-1}q) \longmapsto a^kq^{-\frac{r+k}{2}} \; .
\]
The canonical representation $V$ equals the
standard two-dimensional representation of $GL_{2,\BQ} \,$.
Then,
\[
V_{\ua} = \Sym^k V \bigl( -\frac{r+k}{2} \bigr) \in \Rep (GL_{2,\BQ}) \; ,
\]
where we denote by $\Sym^k V$ the $k$-th symmetric power of $V$.
Theorem~\ref{3A} therefore shows that
\[
{}^{\ua} \CV = \Sym^k \pi_*^1 \one_{\B} \bigl( \frac{-r+k}{2} \bigr)[-r+k] \in \CHQMM \; ,
\]
where $\pi: \B \to M^K$ is the universal elliptic curve. 

The level $K$ equals the kernel of the reduction 
\[
GL_2(\widehat{\BZ}) \longto GL_2 (\BZ / n \BZ) \; ,
\]
for a fixed integer $n \ge 3$.

We claim that the motive $i^* j_* {}^{\ua} \CV \in \DBcM (\partial (M^K)^*)$ is without weights
\[
-(k-1), -(k-2), \ldots, k \; ,
\]
and that both weights $-k$ and $k+1$ \emph{do} occur in $i^* j_* {}^{\ua} \CV$. 
In particular, ${}^{\ua} \CV$ belongs to
$\CHMM_{\BQ,\partial w \ne 0,1} \subset \CHMM_\BQ$
if and only if $k \ge 1$, \emph{i.e.}, if and only if $\ua$ is regular.

In order to show the claim, note first that for fixed $k \ge 0$, its
validity does not
depend on the value of $r \equiv k \mod 2$. We may therefore assume that $r = k$,
\emph{i.e.}, that 
\[
{}^{\ua} \CV = \Sym^k \pi_*^1 \one_{\B} \; .
\]
The Chow motive ${}^{\ua} \CV$ thus equals a direct factor of  
$\pi_{k,*} \one_{\B^k}$, where we denote by 
\[
\pi_k : \B^k  := \B \times_{M^K} \ldots \times_{M^K} \B \longto M^K
\]
the projection of the $k$-fold fibre product $\B^k$ to $M^K$.
Concretely, the symmetric group $\FS_k$ acts on $\B^k$ by permutations,
the $k$-th power of the group $\BZ / n \BZ$ by translations, and the
$k$-th power of the
group $\mu_2$ by inversion in the fibres. Altogether \cite[Sect.~1.1.1]{S}, 
this gives a canonical action of the semi-direct product
\[
\Gamma_k := \bigl( (\BZ / n \BZ)^2 \rtimes \mu_2 \bigr)^k \rtimes \FS_k
\]
by automorphisms on $\B^k$. As in \cite[Sect.~1.1.2]{S}, let $\varepsilon: \Gamma_k \to \{ \pm 1 \}$
be the morphism which is trivial on $(\BZ / n \BZ)^{2k}$, is the
product map on $\mu_2^k$, and is the sign character on $\FS_k$.
Let $e$ denote the idempotent in the group ring
$\BQ[\Gamma_k]$ associated to $\varepsilon \ $:
\[
e := \frac{1}{(2n^2)^k \cdot k!} 
\sum_{\gamma \in \Gamma_k} \varepsilon(\gamma)^{-1} \cdot \gamma 
= \frac{1}{(2n^2)^k \cdot k!} 
\sum_{\gamma \in \Gamma_k} \varepsilon(\gamma) \cdot \gamma  \; .
\]
By passage to the graph,
every endomorphism of the $M^K$-scheme $\B^k$ yields a cycle on $B^k \times_{M^K} B^k$
of dimension $d_{\B^k} = k + 1$. Therefore, we may and do consider $e$ as an 
idempotent of $\ch_{k+1} (\B^k \times_{M^K} \B^k) \otimes_\BZ \BQ \,$. Then,
\[
{}^{\ua} \CV = \Sym^k \pi_*^1 \one_{\B} = \bigl( \pi_{k,*} \one_{\B^k} \bigr)^e \; .
\]
Next, note that $\partial (M^K)^*$, the scheme of cusps of $M^K$, is finite over
$\Spec \BQ \,$. According to Corollary~\ref{2P}, our claim is equivalent to 
the following: the motive $\dMgm(\B^k)^e \in \DBcQQM$
is without weights $- k,- (k - 1),\ldots, k-1$, and both weights $- (k+1)$ and $k$
occur in $\dMgm(\B^k)^e$.

Independently of everything said in this article, this latter claim can be proved
purely geometrically. More precisely, using the detailed analysis of the geometry
of the boundary of the \emph{canonical compactification} of $\B^k$ from \cite[Sect.~2, 3]{S},
one shows that there is an exact triangle 
\[
\Mgm \bigl( \partial (M^K)^* \bigr) (k+1)[k+1] \longto \dMgm(\B^k)^e
\longto \Mgm \bigl( \partial (M^K)^* \bigr)[k] 
\stackrel{[1]}{\longto} 
\]
in $\DBcQQM$ \cite[Ex.~4.16, Rem.~3.5~(b)]{W4}.
But $\Mgm ( \partial (M^K)^* ) (k+1)[k+1]$
is pure of weight $-(k+1)$, and $\Mgm ( \partial (M^K)^* ) [k]$
is pure of weight $k$.

For $k \ge 1$, this shows that $\dMgm(\B^k)^e$
is indeed without weights $- k,- (k - 1),\ldots, k-1$. In particular, 
\[
\Mgm \bigl( \partial (M^K)^* \bigr) (k+1)[k+1] \longto \dMgm(\B^k)^e
\longto \Mgm \bigl( \partial (M^K)^* \bigr)[k] 
\stackrel{[1]}{\longto} 
\]
is \emph{the} weight filtration of $\dMgm(\B^k)^e$ avoiding weights $-1$ and $0$
\cite[Cor.~1.9]{W4}. Since $\Mgm ( \partial (M^K)^* ) \ne 0$, this shows that
both weights $-(k+1)$ and $k$ occur in $\dMgm(\B^k)^e$.

For $k=0$, we have $\dMgm(\B^k)^e = \dMgm(M^K)$. The exact triangle
\[
\Mgm \bigl( \partial (M^K)^* \bigr) (1)[1] \longto \dMgm(M^K)^e
\longto \Mgm \bigl( \partial (M^K)^* \bigr) 
\stackrel{[1]}{\longto} 
\]
is split (cmp.~\cite[Ex.~4.12]{W4}); therefore both weights $-1$ and $0$
occur in $\dMgm(\B^k)^e = \dMgm(M^K)$.
\end{Ex}

\begin{Rem} \label{3D}
(a)~The statements on $\dMgm(\B^k)^e$ from Example~\ref{3C}
admit integral versions. More precisely (see \cite[Sect.~3]{W4}),
they hold in the category $DM_{gm}(\BQ)$ of geometrical
motives over $\BQ \,$,
tensored with $\BZ[1 / (2n \cdot k!)]$. \\[0.1cm]
(b)~The analysis of the geometry
of the boundary of the canonical compacti\-fication of $\B^k$ from \cite[Sect.~2, 3]{S}
can be employed to show directly that the exact triangle
\[
\Mgm \bigl( \partial (M^K)^* \bigr) (k+1)[k+1] \longto \dMgm(\B^k)^e
\longto \Mgm \bigl( \partial (M^K)^* \bigr)[k] 
\stackrel{[1]}{\longto} 
\]
is induced by an exact triangle 
\[
\one_{\partial (M^K)^*}[-k] \longto i^* j_* {}^{\ua} \CV
\longto \one_{\partial (M^K)^*}(-(k+1))[-(k+1)]
\stackrel{[1]}{\longto}
\]
in $\DBcM (\partial (M^K)^*)$. We leave the details to the reader. \\[0.1cm]
(c)~To the author's knowledge,
Example~\ref{3C} is the only non-trivial case of a non-compact pure Shimura variety,
where the weights in $i^* j_* {}^{\ua} \CV$ can be controlled by purely geometrical,
\emph{i.e.}, intrinsically motivic, means. \\[0.1cm]
(d)~For arbitrary level $K$, the motive $i^* j_* {}^{\ua} \CV \in \DBcM (\partial (M^K)^*)$ 
is still without weights
\[
-(k-1), -(k-2), \ldots, k \; ,
\]
and both weights $-k$ and $k+1$ still occur. 
A motivic proof would run as follows. First,
using conjugation in $GL_2(\BA_f)$, one may assume that $K$ contains the kernel
$K_n$ modulo $n$, for an integer $n \ge 3$. Then, one analyzes 
$K / K_n$-equivariance of the direct image of the exact triangle
\[
\one_{\partial (M^{K_n})^*}[-k] \longto i^* j_* {}^{\ua} \CV
\longto \one_{\partial (M^{K_n})^*}(-(k+1))[-(k+1)]
\stackrel{[1]}{\longto}
\]
under the finite morphism from $\partial (M^{K_n})^*$ to $\partial (M^K)^*$.

Alternatively, use the case $L = \BQ$ from Example~\ref{3E}
(which relies on realizations).
\end{Rem}

\begin{Ex} \label{3E}
Our next example concerns \emph{Hilbert--Blumenthal varieties}.
Fix a totally real number field $L$, and denote by $g$ its degree.
The reductive group $G$ equals the fibre product
\[
G := \Gm \times_{\ReL \GLm , \det} \ReL GL_{2,L} 
\]
($\ReL:=$ the Weil restriction from $L$ to $\BQ$).
The reflex field equals $\BQ \,$, and $F = F'$ is a sub-field of $\BR$
containing the images $\sigma(L)$, for $\sigma$ running 
through the set $I_L$ of all real embeddings of $L$.
The dominant character $\ua$ is identified with a $(g+1)$-tuple of integers 
$((k_\sigma)_{\sigma \in I_L},r)$,
with $k_\sigma \ge 0$ and $r \equiv \sum_{\sigma \in I_L} k_\sigma \mod 2$: 
note that 
\[
(\ReL \GL_{2,L}) \times_\BQ F \isoto \prod_{\sigma \in I_L}     
\GL_{2,F} \; ,
\]
and under that isomorphism, 
\[
G \times_\BQ F \isoto 
\biggl\{ (M_\sigma)_{\sigma \in I_L} \in \prod_{\sigma \in I_L} \GL_{2,F} \; , \;
\det(M_\sigma) =  \det(M_\eta) \quad \forall \; \sigma, \eta \in I_L \biggr\} \; .
\]
Choosing $T \subset G \times_\BQ F \subset \prod_{\sigma \in I_L} \GL_{2,F}$ to be equal to 
the sub-group of elements having diagonal entries at each $\sigma \in I_L$, we have
\[
\ua : T \longto \GFm \; , \; 
\bigl( \diag(a_\sigma,a_\sigma^{-1}q) \bigr)_{\sigma \in I_L} \longmapsto 
(\prod_{\sigma \in I_L} a_\sigma^{k_\sigma}) \cdot 
q^{-\frac{r+\sum_\sigma k_\sigma}{2}} \; .
\]
The canonical representation $V$ equals the ($2g$-dimensional) Weil restriction of the
standard two-dimensional representation of $GL_{2,L}$. Thus,
\[
V \times_\BQ F \isoto \bigoplus_{\sigma \in I_L} V_\sigma \; ,
\]
where for $\sigma  \in I_L$, we denote by $V_\sigma$ the standard  
two-dimensional representation of $GL_{2,F}$, seen as the $\sigma$-component
of $(\ReL \GL_{2,L}) \times_\BQ F$. It is the image of an idempotent endomorphism
$e_\sigma$ of $V \times_\BQ F \in \Rep (G_F)$, whose kernel equals
$\oplus_{\eta \ne \sigma} V_\eta$.
Then, setting $k := \sum_\sigma k_\sigma$, 
\[
V_{\ua} = \bigl( \bigotimes_{\sigma \in I_L} \Sym^{k_\sigma} V_\sigma \bigr) 
\bigl( -\frac{r+k}{2} \bigr) \in \Rep (G_F) \; .
\]
Theorem~\ref{3A} therefore shows that
\[
{}^{\ua} \CV =  \bigl( \bigotimes_{\sigma \in I_L} 
\Sym^{k_\sigma} \pi_*^1 \one_{\B}^{e_\sigma} \bigr) 
\bigl( \frac{-r+k}{2} \bigr)
[-r+k] \in \CHFMM \; ,
\]
where $\pi: \B \to M^K$ is the universal Abelian $g$-fold, 
and for every $\sigma \in I_L$, we denote by
$e_\sigma: \pi_*^1 \one_{\B} \to \pi_*^1 \one_{\B}$ the idempotent
endomorphism of $\pi_*^1 \one_{\B}$ induced by functoriality.

We claim that for any (neat) level $K$, the following is true: 
the motive $i^* j_* {}^{\ua} \CV$ is zero unless $\ua$ is \emph{parallel}, \emph{i.e.},
unless all $k_\sigma$ are equal to each other.
Furthermore, $i^* j_* {}^{\ua} \CV$ is
without weights
\[
-(k-1), -(k-2), \ldots, k \; ,
\]
and both weights $-k$ and $k+1$ 
\emph{do} occur in $i^* j_* {}^{\ua} \CV$, provided $\ua$ is parallel. 
In particular, ${}^{\ua} \CV$ belongs to
$\CHMM_{F,\partial w \ne 0,1} \subset \CHMM_F$
if and only if $k = \sum_\sigma k_\sigma \ge 1$, \emph{i.e.}, if and only if
at least one of the tensor components $\Sym^{k_\sigma} V_\sigma$ of $V_{\ua}$ is regular.

In order to show the claim, note first that as in Example~\ref{3C},
we may assume that $r = k$,
\emph{i.e.}, that 
\[
{}^{\ua} \CV = \bigotimes_{\sigma \in I_L} 
\Sym^{k_\sigma} \pi_*^1 \one_{\B}^{e_\sigma} \; .
\]
The Chow motive ${}^{\ua} \CV$ thus equals a direct factor of  
$\pi_{k,*} \one_{\B^k}$, where 
we denote by $\pi_k$
the projection of the $k$-fold fibre product $\B^k$ to $M^K$.
For a concrete description of the associated idempotent endomorphism
\[
e \in \ch_{g(k+1)} (\B^k \times_{M^K} \B^k) \otimes_\BZ F \; ,
\]
we refer to \cite[Lemma~3.4]{W7}. 

Next, by \cite[Thm.~3.5, 3.6]{W7}, the motive $\dMgm(\B^k)^e \in \DBcFQM$
is zero if $\ua$ is not parallel, and it
is without weights $- k,- (k - 1),\ldots, k-1$. 
Now as in Example~\ref{3C}, the scheme of cusps $\partial (M^K)^*$ 
is finite over $\Spec \BQ \,$.
According to Corollary~\ref{2P}, the motive $i^* j_* {}^{\ua} \CV$ is therefore
zero unless $\ua$ is parallel, and it is
without weights $-(k-1), -(k-2), \ldots, k$.

It remains to show that if $\ua$ is parallel, then
both weights $- k$ and $k+1$ occur in $i^* j_* {}^{\ua} \CV$,
or equivalently, both weights $- (k+1)$ and $k$ occur in $\dMgm(\B^k)^e$.
All $k_\sigma$ are equal to each other, say
\[
k_\sigma = s \ge 0 \quad \forall \; \sigma \in I_L \; .
\]
Thus, $k = g \cdot s$, and
\[
V_{\ua} = \bigl( \bigotimes_{\sigma \in I_L} \Sym^s V_\sigma \bigr) ( -k ) \; ,
\]
which is isomorphic to $\bigotimes_{\sigma \in I_L} \Sym^s V_\sigma^\vee$
(where we denote by $V_\sigma^\vee$ the dual of $V_\sigma$, for all $\sigma \in I_L$).

By \cite[Prop.~2.5]{W7}, the motive $\dMgm(\B^k)^e$ realizes
to give \emph{boundary cohomology} 
\[
\bigl( (\partial H^n (\B^k(\BC),\BQ) \otimes_\BQ F)^e \bigr)_{n \in \BZ} \; .
\]
By \cite[Prop.~4.5]{W7}, for all $n \in \BZ$, there are isomorphisms of Hodge structures
\[
( \partial H^n (\B^k(\BC),\BQ) \otimes_\BQ F )^e \isoto
\partial H^{n-k} (M^K(\BC), \mu_{\bf H}(V_{\ua}))
\]
(note that the weight of $V_{\ua}$ equals $k$). 
According to the proof of \cite[Thm.~3.5]{W7}, in particular, \cite[pp.~2351 and 2352]{W7},
\[
\partial H^0 (M^K(\BC), \mu_{\bf H}(V_{\ua}))
\]
is non-zero, and pure of weight $0$, while
\[
\partial H^{2g-1} (M^K(\BC), \mu_{\bf H}(V_{\ua}))
\]
is non-zero, and pure of weight $2(k + g)$. Therefore,
\[
( \partial H^{k} (\B^k(\BC),\BQ) \otimes_\BQ F )^e \ne 0
\]
is pure of weight $0$, and 
\[
( \partial H^{k+2g-1} (\B^k(\BC),\BQ) \otimes_\BQ F )^e \ne 0
\]
is pure of weight $2(k + g)$. Therefore, weights $0-k = -k$ and
$2(k+g) - (k+2g-1) = k+1$ occur in the realization of $\dMgm(\B^k)^e$.
Given that the realization on geometrical motives is contravariant,
and exchanges the signs of weights, this implies in particular that
weights $-(k+1)$ and
$k$ occur in $\dMgm(\B^k)^e$.
\end{Ex}

\begin{Rem} \label{3F}
The proofs of \cite[Thm.~3.5, 3.6]{W7} rely on the fact that $\dMgm(\B^k)$ is a 
\emph{Dirichlet--Tate motive over $\BQ \,$},
and that on such motives, the realizations are weight conservative \cite[Cor.~3.10~(c)]{W8}.
\end{Rem}

\begin{Ex} \label{3G}
The third example concerns \emph{Picard varieties}.
Fix a $CM$-field $L$,  
and denote by $2g$ its degree. 
Fix a three-dimensional $L$-vector space $V'$, together with an $L$-valued
non-degenerate Hermitian form $J$, such that for every 
$\sigma$ in the set $I_L$ of complex embeddings of $L$, the form
$J_\sigma := J \otimes_{L,\sigma} \BC$ is of signature $(2,1)$.
Fix a $CM$-type $\Psi$ of $L$; thus, the set $I_L$ is the disjoint union of $\Psi$ 
and of its conjugate.  
The reductive group $G$ equals the group of unitary similitudes
\[
G := GU(V',J) \subset \ReL \GL_L(V') \; .
\]
Thus, for any $\BQ$-algebra $R$, the group $G(R)$ equals
\[
\biggl\{ g \in \GL_{L \otimes_\BQ R} (V' \otimes_\BQ R) \; , \; 
\exists \, \lambda(g) \in R^* \; , \; 
J(g \argdot,g \argdot) = \lambda(g) \cdot J(\argdot,\argdot) \biggr\} \; .
\] 
In particular, the similitude norm $\lambda(g)$ defines a canonical morphism
\[
\lambda: G \longto \Gm \; .
\]
Then $F = F'$ is a sub-field of $\BC$
containing the images $\sigma(L)$, for all $\sigma \in I_L$.
There is an isomorphism
\[
G \times_\BQ F \isoto \bigl( \prod_{\sigma \in \Psi} \GL_F (V_F') \bigr) \times_F \GFm
\]
(\cite[p.~4]{C}; see \cite[p.~364]{W7a} for the case $g=1$).

Fix a basis of $V'$; thus,
\[
G \times_\BQ F \isoto \bigl( \prod_{\sigma \in \Psi} \GL_{3,F} \bigr) \times_F \GFm \; .
\]
The dominant character $\ua$ is identified with a 
$(3g+1)$-tuple of integers 
$((k_{1,\sigma},k_{2,\sigma},c_\sigma)_{\sigma \in \Psi},r)$,
with $k_{1,\sigma} \ge k_{2,\sigma} \ge 0$,
\[
c_\sigma \equiv k_{1,\sigma} + k_{2,\sigma} \mod 2 \; , \quad \text{and} \quad 
r \equiv \sum_{\sigma \in \Psi} \frac{c_\sigma + k_{1,\sigma} + k_{2,\sigma}}{2} \mod 2 \; :
\] 
choosing $T \subset ( \prod_{\sigma \in \Psi} \GL_{3,F} ) \times_F \GFm$ to be equal to 
the sub-group of elements having diagonal entries at each $\sigma \in \Psi$, 
the character $\ua$ sends 
\[
\bigl( \bigl( \diag(a_\sigma,a_\sigma^{-1}b_\sigma,b_\sigma^{-1} q_\sigma) 
\bigr)_{\sigma \in \Psi} ,
f \bigr) \in T(F)
\]
to
\[
\bigl( \prod_{\sigma \in \Psi} a_\sigma^{k_{1,\sigma} - k_{2,\sigma}} b_\sigma^{k_{2,\sigma}}
q_\sigma^{\frac{c_\sigma - (k_{1,\sigma} + k_{2,\sigma})}{2}} \bigr) \cdot 
f^{-\halb (r + \sum_\sigma \frac{3c_\sigma - (k_{1,\sigma} + k_{2,\sigma})}{2})}  
\in \BG_m(F) \; .
\]
Under this normalization, the similitude norm $\lambda$ corresponds to the character
$\ua = ((0,0,0)_{\sigma \in \Psi},-2)$. The corresponding irreducible representation
$V_\lambda$ of $G_F$, and actually, of $G$, is the Tate twist $\BQ(1)$. 
It follows from Theorem~\ref{3A}~(c) that as far as control of weights is concerned,
we may replace a given $\ua = ((k_{1,\sigma},k_{2,\sigma},c_\sigma)_{\sigma \in \Psi},r)$
by $((k_{1,\sigma},k_{2,\sigma},c_\sigma)_{\sigma \in \Psi},r')$, with 
$r' \equiv r \mod 2$.  

\forget{
The canonical representation $V$ equals the ($6g$-dimensional) Weil restriction of $V'$. 
Thus,
\[
V \times_\BQ F \isoto \bigoplus_{\sigma \in I_L} V_\sigma \; ,
\]
where for $\sigma  \in I_L$, we denote by $V_\sigma$ the standard  
two-dimensional representation of $GL_{2,F}$, seen as the $\sigma$-component
of $(\ReL \GL_{2,L}) \times_\BQ F$. It is the image of an idempotent endomorphism
$e_\sigma$ of $V \times_\BQ F \in \Rep (G_F)$, whose kernel equals
$\oplus_{\eta \ne \sigma} V_\eta$.
Then, setting $k := \sum_\sigma k_\sigma$, 
\[
V_{\ua} = \bigl( \bigotimes_{\sigma \in I_L} \Sym^{k_\sigma} V_\sigma \bigr) 
\bigl( -\frac{r+k}{2} \bigr) \in \Rep (G_F) \; .
\]
Theorem~\ref{3A} therefore shows that
\[
{}^{\ua} \CV =  \bigl( \bigotimes_{\sigma \in I_L} 
\Sym^{k_\sigma} \pi_*^1 \one_{\B}^{e_\sigma} \bigr) 
\bigl( \frac{-r+k}{2} \bigr)
[-r+k] \in \CHFMM \; ,
\]
where $\pi: \B \to M^K$ is the universal Abelian $g$-fold, 
and for every $\sigma \in I_L$, we denote by
$e_\sigma: \pi_*^1 \one_{\B} \to \pi_*^1 \one_{\B}$ the idempotent
endomorphism of $\pi_*^1 \one_{\B}$ induced by functoriality.
}

We claim that  
the motive ${}^{\ua} \CV$ belongs to
$\CHMM_{F,\partial w \ne 0,1}$
if $\ua$ is regular, \emph{i.e.}, if for all $\sigma \in \Psi$, we have
\[
k_{1,\sigma} > k_{2,\sigma} > 0 \; .
\]
In order to show the claim, 
it will be useful to compare our parametrization of characters to that
of \cite[Sect.~3]{C}. There, the standard basis of characters of
the split torus $T$ is used. First, write 
$((a_\sigma,b_\sigma,\gamma_\sigma)_{\sigma \in \Psi},d)$
instead of $((a_\sigma,b_\sigma,c_\sigma)_{\sigma \in \Psi},d)$ as in \loccit ,
and $\lambda((a_\sigma,b_\sigma,\gamma_\sigma)_{\sigma \in \Psi},d)$ 
for the associated character. We then leave it to the reader to
verify that the change of parameters 
$\lambda((a_\sigma,b_\sigma,\gamma_\sigma)_{\sigma \in \Psi},d) 
\leftrightarrow ((k_{1,\sigma},k_{2,\sigma},c_\sigma)_{\sigma \in \Psi},r)$
is the following: the character
$((k_{1,\sigma},k_{2,\sigma},c_\sigma)_{\sigma \in \Psi},r)$ of $T$ equals
\[ 
\lambda \bigl( 
\bigl( \frac{c_\sigma+k_{1,\sigma}-k_{2,\sigma}}{2},
\frac{c_\sigma-k_{1,\sigma}+k_{2,\sigma}}{2},
\frac{c_\sigma-(k_{1,\sigma}+k_{2,\sigma})}{2} \bigr)_{\sigma \in \Psi}, d
\bigr) \; ,
\]
with $d = -\halb (r + \sum_\sigma \frac{3c_\sigma-(k_{1,\sigma}+k_{2,\sigma})}{2})$.
The character $\lambda((a_\sigma,b_\sigma,\gamma_\sigma)_{\sigma \in \Psi},d)$ equals
\[
\bigl( \bigl( a_\sigma-\gamma_\sigma,b_\sigma-\gamma_\sigma,
                         a_\sigma+b_\sigma \bigr)_{\sigma \in \Psi},
-(2d+\sum_\sigma(a_\sigma+b_\sigma+\gamma_\sigma)) \bigr) \; .
\]
The character $\lambda((a_\sigma,b_\sigma,\gamma_\sigma)_{\sigma \in \Psi},d)$
is dominant if and only if for for all $\sigma \in \Psi$,
\[
a_\sigma \ge b_\sigma \ge \gamma_\sigma \; ,
\]
and it is regular if and only if for for all $\sigma \in \Psi$,
\[
a_\sigma > b_\sigma > \gamma_\sigma \; .
\]
To replace $\ua = ((k_{1,\sigma},k_{2,\sigma},c_\sigma)_{\sigma \in \Psi},r)$
by $((k_{1,\sigma},k_{2,\sigma},c_\sigma)_{\sigma \in \Psi},r')$, with 
$r' \equiv r \mod 2$, means to replace 
$\lambda((a_\sigma,b_\sigma,\gamma_\sigma)_{\sigma \in \Psi},d)$
by $\lambda((a_\sigma,b_\sigma,\gamma_\sigma)_{\sigma \in \Psi},d')$. 
This together with \cite[Prop.3.2]{C} shows that 
as far as our weight estimates are concerned,
we may assume, by modifying the value of $r$ if necessary, 
that the Chow motive ${}^{\ua} \CV$ equals a direct factor of  
$\pi_{r,*} \one_{\B^r}$, where 
we denote by $\pi_r$
the projection of the $r$-fold fibre product $\B^r$ to $M^K$.
Let $e$ be the associated idempotent endomorphism of $\pi_{r,*} \one_{\B^r}$. 

By \cite[Thm.~3.6]{C}, the motive $\dMgm(\B^r)^e \in \DBcFEM$
is without weights $- 1$ and $0$ if $\ua$ is regular. 
Now as in Examples~\ref{3C} and \ref{3E}, the scheme of cusps $\partial (M^K)^*$ 
is finite over $\Spec E$.
According to Corollary~\ref{2P}, the motive 
${}^{\ua} \CV$ therefore lies in $\CHMM_{F,\partial w \ne 0,1}$
if $\ua$ is regular.
\end{Ex}

\begin{Rem} \label{3H}
(a)~For $g=1$, it is shown in \cite[Thm.~3.8]{W7a} that regu\-larity of $\ua$
is not only sufficient, but also necessary for ${}^{\ua} \CV$
to belong to $\CHMM_{F,\partial w \ne 0,1}$. 
In fact, a precise interval of weights avoided
by $i^* j_* {}^{\ua} \CV$ is given: putting 
$k := \min(k_{1,\sigma}-k_{2,\sigma},k_{2,\sigma})$ (for the unique $\sigma$ in $\Psi$),
the motive $i^* j_* {}^{\ua} \CV$ is without weights
\[
-(k-1), -(k-2), \ldots, k \; ,
\]
and both weights $-k$ and $k+1$ 
do occur in $i^* j_* {}^{\ua} \CV$
(\cite[Thm.~3.8]{W7a}, Corollary~\ref{2P}). \\[0.1cm]
(b)~As soon as $g \ge 2$, regularity of $\ua$
is no longer necessary for ${}^{\ua} \CV$
to belong to $\CHMM_{F,\partial w \ne 0,1}$,
as illustrated by \cite[Prop.~3.8]{C}. 
By checking the details of the computations leading to the proof of
\cite[Thm.~3.6]{C} (see in particular \cite[pp.~21--24]{C}), one sees that 
the motive ${}^{\ua} \CV$
belongs to $\CHMM_{F,\partial w \ne 0,1}$ if and only if 
at least one of the components $(k_{1,\sigma},k_{2,\sigma},c_\sigma)$
of $\ua$ is regular, \emph{i.e.}, if and only if 
\[
k_{1,\sigma} > k_{2,\sigma} > 0 
\]
for \emph{some} $\sigma \in \Psi$.
Note that this condition is equivalent to saying that $\ua$ is not
\emph{Kostant-parallel for $m=-1$} in the terminology of \cite[Def.~3.2]{C}.

It is proved \cite[Prop.~3.7]{C} that $i^* j_* {}^{\ua} \CV \ne 0$ if and only if
$\ua$ is Kostant-parallel for some integer $m \ge -1$
(this is the case in particular if $\ua$ is parallel,
\emph{i.e.}, if all pairs $(k_{1,\sigma},k_{2,\sigma})$ are equal to each other).
Thus, assume $i^* j_* {}^{\ua} \CV \ne 0$, and denote by $m$ the minimal value
for which $\ua$ is Kostant-parallel.
Define $n$ as the number of $\sigma \in \Psi$ for which 
$k_{1,\sigma} > m$.
A close analysis of \cite[pp.~21--24]{C} reveals that 
$i^* j_* {}^{\ua} \CV$ is without weights
\[
-(mg+n-1), -(mg+n-2), \ldots, mg+n \; ,
\]
and that both weights $-(mg+n)$ and $mg+n+1$ 
do occur in $i^* j_* {}^{\ua} \CV$. \\[0.1cm]
(c)~The proofs of \cite[Thm.~3.8]{W7a} and \cite[Thm.~3.6]{C} 
rely on the fact that $\dMgm(\B^r)$ is a motive of Abelian type over $E$,
and that on such motives, the realizations are weight conservative \cite[Thm.~1.13]{W7a}.
\end{Rem}

\begin{Rem} 
None of the examples treated so far necessitates the use of 
the new criteria on absence of weights in the boundary proved in Section~\ref{2}.
Indeed, weights were controlled, respectively, by purely geometrical means (Example~\ref{3C}),
by using weight conservativity of the restriction of the realizations to Dirichlet--Tate
motives (Example~\ref{3E}), and by using weight conservativity of the restriction of the realizations to motives of Abelian type over a point (Example~\ref{3G}). 

All of them \emph{could} be proved using Corollary~\ref{2Ca}, since the main 
assumption: $i^* j_* {}^{\ua} \CV$ of Abelian type, is satisfied.
\end{Rem}

\begin{Ex} \label{3I}
Our last example concerns \emph{Hilbert--Siegel varieties}.
Fix a totally real number field $L$, and denote by $g$ its degree.
Fix a four-dimensional $L$-vector space $V'$, together with an $L$-valued
non-degenerate symplectic bilinear form $J$. 
The reductive group $G$ equals the group of symplectic similitudes
\[
G := GSp(V',J) \subset \ReL \GL_L(V') \; .
\]
Thus, for any $\BQ$-algebra $R$, the group $G(R)$ equals
\[
\biggl\{ g \in \GL_{L \otimes_\BQ R} (V' \otimes_\BQ R) \; , \; 
\exists \, \lambda(g) \in R^* \; , \; 
J(g \argdot,g \argdot) = \lambda(g) \cdot J(\argdot,\argdot) \biggr\} \; .
\] 
In particular, the similitude norm $\lambda(g)$ defines a canonical morphism
\[
\lambda: G \longto \Gm \; .
\]
The reflex field equals $\BQ \,$, and $F = F'$ is a sub-field of $\BR$
containing the images $\sigma(L)$, for $\sigma$ running 
through the set $I_L$ of all real embeddings of $L$.

Fix a symplectic basis of $V$, 
in which $J$ acquires the $4 \times 4$-matrix
\[
\left( \begin{array}{cc}
0 & I_2 \\
-I_2 & 0
\end{array} \right) \; ;
\]  
thus,
\[
G \times_\BQ F \longinto \ReL GSp_{4,L} \; ,
\]
where $GSp_{4,L}$ denotes the sub-group of 
$GL_{4,\BQ}$ of matrices $g$
satisfying the relation
\[
 ^{t}\! g J g = \lambda(g) \cdot J \; .
\]
More precisely, 
\[
G := \Gm \times_{\ReL \GLm , \det} \ReL GSp_{4,L} \; . 
\]
The dominant character $\ua$ is identified with a $(2g+1)$-tuple of integers 
$((k_{1,\sigma},k_{2,\sigma})_{\sigma \in I_L},r)$,
with $k_{1,\sigma} \ge k_{2,\sigma} \ge 0$
and $r \equiv \sum_{\sigma \in I_L} (k_{1,\sigma} + k_{2,\sigma}) \mod 2$: 
choosing $T \subset G \times_\BQ F \subset \prod_{\sigma \in I_L} GSp_{4,F}$ to be equal to 
the sub-group of elements having diagonal entries at each $\sigma \in I_L$, 
the character $\ua$ sends 
\[
\bigl( \diag(a_\sigma,b_\sigma,a_\sigma^{-1}q,b_\sigma^{-1}q) \bigr)_{\sigma \in I_L} \in T(F) 
\]
to
\[
\bigl( \prod_{\sigma \in I_L} a_\sigma^{k_{1,\sigma}} \bigr) \cdot
\bigl( \prod_{\sigma \in I_L} b_\sigma^{k_{2,\sigma}} \bigr) \cdot
q^{-\frac{r+\sum_\sigma (k_{1,\sigma} + k_{2,\sigma})}{2}} \; .
\]

Using the criterion from Corollary~\ref {2Ca}~(b),
it is shown in \cite[Cor.~2.1.0.4]{Ca} that
the motive ${}^{\ua} \CV$ belongs to
$\CHMM_{F,\partial w \ne 0,1}$
if at least one of the components $(k_{1,\sigma},k_{2,\sigma})$
of $\ua$ is regular, \emph{i.e.}, if  
\[
k_{1,\sigma} > k_{2,\sigma} > 0 
\]
for \emph{some} $\sigma \in I_L$.
\end{Ex}

\begin{Rem}
(a)~For $g=1$, it is shown in \cite[Thm.~1.6]{W12} that regu\-larity of $\ua$
is not only sufficient, but also necessary for ${}^{\ua} \CV$
to belong to $\CHMM_{F,\partial w \ne 0,1}$. 
In fact, a precise interval of weights avoided
by $i^* j_* {}^{\ua} \CV$ is given: putting 
$k := \min(k_{1,\sigma}-k_{2,\sigma},k_{2,\sigma})$ (for the unique $\sigma$ in $I_L$),
the motive $i^* j_* {}^{\ua} \CV$ is without weights
\[
-(k-1), -(k-2), \ldots, k \; ,
\]
and both weights $-k$ and $k+1$ 
do occur in $i^* j_* {}^{\ua} \CV$. \\[0.1cm]
(b)~As soon as $g \ge 2$, regularity of one of the components of $\ua$
is no longer necessary for ${}^{\ua} \CV$
to belong to $\CHMM_{F,\partial w \ne 0,1}$.
Indeed, according to \cite[Cor.~2.1.0.4]{Ca},
for the motive ${}^{\ua} \CV$
not to belong to $\CHMM_{F,\partial w \ne 0,1}$, it is ne\-cessary
and sufficient that none of the components $(k_{1,\sigma},k_{2,\sigma})$
of $\ua$ is regular, and that in addition all $k_{2,\sigma}$ are equal to each other.

For $g = 2$ one might thus choose $\ua$ equal to
$((1,0),(1,1),3)$, to obtain an example of a character none of whose components
is regular, but whose associated motive ${}^{\ua} \CV$
belongs nonetheless to $\CHMM_{F,\partial w \ne 0,1}$. \\[0.1cm]
(c)~In \cite[Thm.~2.1.0.3]{Ca}, an interval of weights avoided
by $i^* j_* {}^{\ua} \CV$, sharp in most cases, 
is given for any value of $g$.
As in the case of Picard varieties, the notion of Kostant-parallelism occurs.
But in addition, the notion of \emph{corank} is needed.
We refer to \loccit \ for the definition of the corank, and for the precise formulae. \\[0.1cm]
(d)~As was the case in our earlier examples, the above weight estimates are valid 
for \emph{any} neat level $K$.
\end{Rem}

\begin{Rem}
According to Corollary~\ref{2Ca}~(c) and Proposition~\ref{2K},
the weight estimates for $i^* j_* {}^{\ua} \CV$ induce weight estimates
for the corresponding $e$-part $\dMgm(\B^r)^e$ of the boundary motive $\dMgm(\B^r)$. 

Contrary to the cases treated in Examples~\ref{3C}, \ref{3E} and \ref{3G},
Hilbert--Siegel varieties have a boundary $\partial (M^K)^*$ of strictly positive
dimension (equal to the degree $g$ of the totally real number field $L$).
More seriously, that boundary is not of Abelian type over $\Spec \BQ \,$.
Thus, our earlier results on weight conservativity \cite{W8,W7a} cannot be employed
to control the weights in $\partial (M^K)^*$ directly from those in the 
boundary cohomology of $M^K$.

The strategy set out in Section~\ref{2}, \emph{i.e.}, the detour
\emph{via} relative motives over $(M^K)^*$, \emph{is} therefore needed 
in order to treat Example~\ref{3I}.
\end{Rem}

Given our examples, the following seems to be justified.

\begin{Ques}[{\cite[Conj.~A]{W13}}] \label{3J}
Let $M^K$ be a Shimura variety of $PEL$-type, associated to pure Shimura
data $(G,\CH)$ and a neat level $K$. Denote by $j$
the open immersion of $M^K$ into its Baily--Borel 
compactification $(M^K)^*$, and by $i$ the closed
immersion of the boun\-dary $\partial (M^K)^*$ of $(M^K)^*$. 
Let $V_{\ua}$ be an irreducible regular representation of $G$. \\

Does the motive ${}^{\ua} \CV = \widetilde{\mu}(V_{\ua})$ belong to
$\CHMM_{F,\partial w \ne 0,1} \,$?
\end{Ques}

Our examples actually suggest that weaker conditions on $\ua$ are still 
sufficient for ${}^{\ua} \CV$ to belong to
$\CHMM_{F,\partial w \ne 0,1}$. We refer to \cite{W13} for stronger versions
of Question~\ref{3J}. 


\bigskip

%
%

\end{document}